\documentclass[letterpaper,10pt,conference]{ieeeconf}

\IEEEoverridecommandlockouts  
\usepackage{graphicx}
\usepackage{dsfont} 
\usepackage{amsmath}
\usepackage{amssymb}
\usepackage{url}
\usepackage{algpseudocode}
\usepackage{color}
\usepackage[utf8]{inputenc} 
\usepackage[format=hang,singlelinecheck=false,caption=false,font=footnotesize]{subfig}
\usepackage{trimclip}

\graphicspath{{./figures/}}

\let\originalleft\left
\let\originalright\right
\DeclareRobustCommand{\left}{\mathopen{}\mathclose\bgroup\originalleft}
\DeclareRobustCommand{\right}{\aftergroup\egroup\originalright}

\newcommand{\R}{\mathbb{R}}

\renewcommand{\epsilon}{\varepsilon}
\renewcommand{\phi}{\varphi}
\newcommand{\derivoper}{\mathrm{d}}

\newcommand{\abs}[1]{\left\lvert #1\right\rvert}

\newcommand{\norm}[1]{\left\lVert #1\right\rVert}

\newcommand{\set}[1]{\left\{ #1\right\}}

\newcommand{\paren}[1]{\left( #1\right)}
\newcommand{\parenb}[1]{\bigl( #1\bigr)}

\newcommand{\sqparen}[1]{\left[ #1\right]}

\newcommand{\diff}[1]{\derivoper #1}

\newcommand{\sqmatrix}[1]{%
  \begin{bmatrix}#1 \end{bmatrix}%
}
\newcommand{\sqmatrixs}[1]{%
  \sqparen{\begin{smallmatrix}#1 \end{smallmatrix}}%
}

\newcommand{\axis}[1]{\hat{#1}}
\newcommand{\obs}[1]{\tilde{#1}}

\newcommand{\topfig}[1]{\fbox{\includegraphics[width=.9in,trim=160 160 60 60,clip]{#1}}}
\newcommand{\botfig}[1]{\fbox{\includegraphics[width=.9in,trim=130 100 6 3,clip]{#1}}}

\begin{document}
\title{%
  Simultaneous Receding Horizon Estimation and Control of a Fencing Robot using a Single Camera%
}

\author{%
  Ignacio de Erausquin \and
  Humberto Gonzalez%
  \thanks{%
    The authors are with the Department of Electrical \& Systems Engineering, Washington University in St.\ Louis, St.\ Louis, MO 63130.
    Emails: {\scriptsize \texttt{ignacio.de.erausquin@gmail.com}, \texttt{hgonzale@wustl.edu}}
  }
}

\thispagestyle{empty}
\pagestyle{empty}
\maketitle

\begin{abstract}
  We present a method for simultaneous Receding Horizon Estimation and Control of a robotic arm equipped with a sword in an adversarial situation.
  Using a single camera mounted on the arm, we solve the problem of blocking a opponent's sword with the robot's sword.
  Our algorithm uses model-based sensing to estimate the opponent's intentions from the camera's observations, while it simultaneously applies a control action to both block the opponent's sword and improve future camera observations.
\end{abstract}

\section{Introduction}
\label{sec:intro}

Today's powerful sensors, such as integrated inertial-measurement units and high-speed cameras, together with efficient embedded processors, allow robots to perceive the environment with a level of detail that is at least comparable to human perception capabilities.
But even using precise and fast sensors, robots are still several steps behind the capabilities of humans at solving everyday tasks.
Indeed, even though prototype robots performing tasks such as driving a car~\cite{Montemerlo2008,Sprinkle2009} or doing chores~\cite{Maitin-Shepard2010} exist today, their performance is still not comparable to that of a human when provided with the same perceptual stimulus.
Among the reasons why robots cannot perform similarly to humans when confronted with similar perceptual stimuli, two arguments stand above the rest:
first, humans perform \emph{active perception}~\cite{Bajcsy1988}, i.e., we actively seek to improve the quality of the data received by our sensors, for example, when we move our head to observe a wider field-of-view;
second, humans can extract the \emph{actionable information}~\cite{Soatto2009,Soatto2011} from continuous streams of data in a natural way, i.e., a loud noise (continuous sound signal) is quickly transformed into a clue of possibility of danger (discrete symbol) by our brain~\cite{Bajcsy1995}.

In this paper we design a strategy to allow a robotic arm equipped with a sword to block its opponents attacks~\cite{fencingrules}.
Using a single camera mounted on one of the links of the robotic arm as the only perceptive information, we provide a practical answer to the problems of performing \emph{active perception} and extracting the \emph{actionable information} from the captured data.
We chose to use a single camera as a sensor since it is passive, small, easy to calibrate, and it is theoretically possible to reconstruct 3-D information using the motion of the arm~\cite{Ma2004}.

We implement two methods, one based on Receding Horizon Control (RHC) to move the robotic arm as a function of the position of the opponent's sword, and one based on Receding Horizon Estimation (RHE) to obtain the 3-D coordinates of the opponent's sword given the images captured by the camera.
More importantly, the RHC method performs \emph{active perception} by adjusting the position of the arm to improve the point of view of the camera, and our RHE method extracts the \emph{actionable information} by using a dynamical model of the opponent's movements.

Several groups have studied the interactions between robots and humans, and have proposed ways to extract information from motion signals~\cite{Goodrich2007,Hanai2011,Tenorth2012,Tenorth2012a}.
In the context of their papers, our results should be understood as a particular approach towards the information retrieval using dynamical models and optimization-based algorithms.

Several groups have used Receding Horizon Control \& Estimation to process data from sensors and solve problems in robotics.
RHE methods have been used extensively in the literature as state observers~\cite{Muske1993,Alessandri2003,Alessandri2010}.
Our approach, as described in Section~\ref{sec:observer}, is not to observe all the opponent's state variables, but rather to estimate the opponent's speed and direction of attack.
Even though the problems are similar, our problem involves a smaller number of variables, greatly simplifying the numerical computations.
Prazenica et al.~\cite{Prazenica2006} use several cameras to estimate the position of moving robots controlled by RHC.
Instead of taking that approach, which would imply that the robotic arm is equipped with a collection of cameras all around the opponent, we preferred the approach where a single camera is actively moved, allowing for a compact, more realistic scenario.

Chipalkatty et al.~\cite{Chipalkatty2010,Chipalkatty2013} consider the problem of incorporating human inputs into the control loop.
Even though in our problem the opponent is assumed to be a human subject, we focus on the design of a completely autonomous fencing robot.

Kunz et al.~\cite{Kunz2011} design a full hybrid dynamical model describing the fencing game, where discrete modes represent the different intentions by the opponents, and the continuous variables describe the configuration space of the two arms.
Their result is complementary to our paper.
While they focus on a complete description of the fencing game without considering how the opponents observe each other, we focus on the interplay between perception and actions, but using only a defensive strategy for the robotic arm.

Liang et al.~\cite{Liang2006} also designed a robot for fencing training.
However, their result focuses on the components of the problem that do not involve moving the robot's sword, but rather emphasize training the fencer's footwork and speed.
Our paper deals entirely with the problem of controlling the sword, thus we do not consider the lateral movement of the two fencers.

It is also worth noting that we currently do not claim any theoretical results regarding the observability of the opponent's sword using a single camera, as studied by Hernandez and Soatto~\cite{Hernandez2013}.
Instead, we try to focus on the proof-of-concept where information from a camera can be naturally transformed into useful (or actionable) information using RHE.

Our paper is organized as follows.
Section~\ref{sec:probdesc} presents the problem of designing a fencing robot in detail, as well as our simplifying assumptions.
Section~\ref{sec:observer} shows the implementation of the RHE method which receives the camera-frame position of the opponent's sword and outputs the speed and direction of the opponent's attack.
Section~\ref{sec:controller} shows the implementation of the RHC method which moves the arm to both defend and improve the point of view of the camera.
Section~\ref{sec:examples} shows simulations of the closed-loop system and analyzes its performance.
Finally, Section~\ref{sec:conclusion} presents our conclusions and future work directions.

\section{Problem Description}
\label{sec:probdesc}

In this section, we describe the model for sword fighting, or \emph{fencing}, between two players, called \emph{fencers}.
We further describe the mathematical model of our robotic arm and camera and how these are used in the fencing problem.

Our model for fencing is based on a simplification of the Olympic saber fencing rules~\cite{fencingrules}.
One fencer, the \emph{attacker}, attempts to perform an attack, while the other fencer, the \emph{defender}, attempts to interrupt this attack.
An attack is described as an uninterrupted motion of the attacker's sword towards the target area.
An attack is considered to be interrupted if the defender forces contact between the two fencers' swords, or if the attacker stops his forward motion towards the target.
In this paper, we will only consider the problem from the defender's perspective. Therefore, the robot will always be the defender, while the opponent will be the attacker.

\subsection{Robot Model}

\begin{figure}[tp]
  \centering
  \resizebox{.9\linewidth}{!}{%
    \input{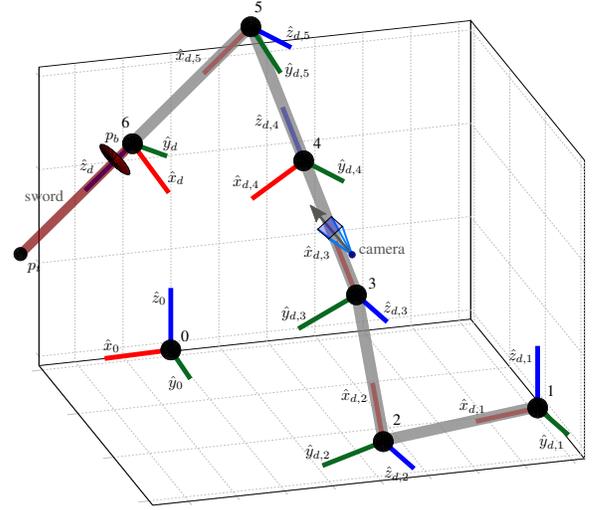}%
  }
  \caption{%
    Diagram of the robotic arm, the coordinate frames defined by each of the joints, the global coordinate frame, the location of the camera, and the location of the sword.
    Each joint rotates around the $\axis{z}$-axis (in color blue).
  }
  \label{fig:robotic_arm}
\end{figure}

The fencing robot is an arm consisting of 6 rotary joints, with corresponding coordinate frames denoted $\paren{\axis{x}_{d,i}, \axis{y}_{d,i}, \axis{z}_{d,i}} \in \R^3$ for each $i \in \set{1,\dotsc,6}$, as shown in Figure~\ref{fig:robotic_arm}.
The subscript $d$ indicates that these frames correspond to the \emph{defender}.
To simplify our notation we will denote the sixth-joint coordinate frame simply as $\paren{\axis{x}_{d}, \axis{y}_{d}, \axis{z}_{d}}$.
The coordinate frames are set such that the \mbox{$i$-th} joint rotates around $\axis{z}_{d,i}$ with angle $\theta_{d,i} \in \R$.
Also, let $\paren{\axis{x}_0, \axis{y}_0, \axis{z}_0}$ be the \emph{world frame}, which is assumed static.
All points in the paper will be represented in the world coordinate frame unless otherwise noted.

We describe the displacement of points in 3-D space using the space of special Euclidean 3-D transformations, $SE(3)$, defined by:
\begin{equation}
	SE(3) = \set{g = (R,T) \in SO(3) \times \R^3}, 
\end{equation}
where $T \in \R^3$ represents a translation and $R \in SO(3)$ represents a rotation, typically described using a \emph{rotation matrix} which leads to the following definition:
\begin{equation}
	SO(3) = \set{R \in \R^{3 \times 3} \mid R^T\, R = I,\, \det(R) = +1}. 
\end{equation}
In the case of our robot, we define 6 transformations, denoted $\set{g_{d,i}}_{i=1}^6 \subset SE(3)$, each transforming points from the \mbox{$(i-1)$-th} frame to the \mbox{$i$-th} frame, where we abuse notation and set \mbox{$0$-th} frame as the world frame.
Note that, since the robot's segments have a fixed length, each of these transformations is parameterized by the joint's angle, i.e., $g_{d,i}(\theta_{d,i})$ for each $i \in \set{1,\dotsc,6}$.

Throughout the paper we use the \emph{homogeneous representation} for points in space and elements in $SE(3)$, as described in Chapter 2 of~\cite{Ma2004}.
Given a point $p \in \R^3$ and a transformation $g = (R,T) \in SE(3)$, then when we apply $g$ to $p$ we get a new point $q = R\, p + T$.
The homogeneous representation summarizes the algebraic operation above by embedding $p$ into $\R^4$ as $\sqmatrixs{p\\1}$ and embedding $g$ into $\R^{4 \times 4}$ as $\sqmatrixs{R & T\\ 0 & 1}$, thus:
\begin{equation}
  \label{eq:hom_rep}
  \sqmatrix{q\\1} = \sqmatrix{R & T\\ 0 & 1}\, \sqmatrix{p\\1}.
\end{equation}
We will abuse notation and assume that all points in $\R^3$ and transformations in $SE(3)$ are in fact written in their homogeneous representation, e.g., we will write equation~\eqref{eq:hom_rep} simply as $q = g\, p$.	

The defender's sword extends along the $\axis{z}_d$-axis, as shown in Figure~\ref{fig:robotic_arm}, which in world coordinates is given by:
\begin{equation}
  \axis{z}_{d}\paren{\Theta_d} = G_d(\Theta_{d})\, \sqmatrixs{0\\0\\1},
\end{equation}
where:
\begin{equation}
	G_d(\Theta_d) = \prod_{i=1}^6 g_{d,i}(\theta_{d,i}),
\end{equation}
and $\Theta_d = \sqmatrix{\theta_{d,1}, \dotsc, \theta_{d,6}}^T$.
Similarly, the base and tip of the defender's sword are described in world coordinates by $G_d(\Theta_d)\, p_b$ and $G_d(\theta_d)\, p_t$, respectively, where:
\begin{equation}
  \label{eq:pt_pb}
  p_b = \sqmatrixs{0\\0\\0},\quad \text{and} \quad
  p_t = \sqmatrixs{0\\0\\\ell}
\end{equation}
are represented in sixth-joint coordinates, and $\ell$ is the length of the sword.

We model the robot's target area, i.e., the area being defended from the attacker, as a finite set of points of interest denoted $S_{\text{target}}$.

We use a kinematic model for the robotic arm, and we assume that the arm is fully actuated.
Hence, the joint angles satisfy the following ODE: $\dot{\Theta}_d(t) = u_d(t)$, where $u_d(t) \in \R^6$ are the angular velocities of the joints.

\subsection{Camera Model}

The fencing robot observes its environment through a single camera mounted on the third joint of the arm, as shown in Figure~\ref{fig:robotic_arm}.
We assume that the camera is capable of obtaining the 2-D pixel coordinates of the tip and base of the attacker's sword.
Thus, in this paper we do not consider the segmentation problem, but it is worth noting that such problem can be solved, for example, by using appropriate color markers.

Given a point $p \in \R^3$ represented in a reference frame $\paren{\axis{x},\axis{y},\axis{z}}$, we will denote each coordinate of $p$ by $\sqparen{p}_{\axis{x}}$, $\sqparen{p}_{\axis{y}}$, and $\sqparen{p}_{\axis{z}}$, respectively.

Let $g_{\text{cam}} \in SE(3)$ denote the transformation from the camera coordinate frame, $\paren{\axis{x}_{\text{cam}},\axis{y}_{\text{cam}},\axis{z}_{\text{cam}}}$, to the third-joint coordinate frame.
The camera coordinate frame is assumed to be located at the center of the camera focal plane.
Hence, the camera's transformation relative to the world frame is:
\begin{equation}
	G_{\text{cam}}(\Theta_d) =  g_{d,1}(\theta_{d,1})\, g_{d,2}(\theta_{d,2})\, g_{d,3}(\theta_{d,3})\, g_{\text{cam}}.
\end{equation}
The camera coordinate frame is only 2-D, hence the observation $\obs{p} \in \R^2$ of a point $p \in \R^3$ is modeled as:
\begin{equation}
  \sqparen{p}_{\axis{z}_{\text{cam}}}\, \obs{p} = K\, \Pi_0\, G_{\text{cam}}^{-1}(\Theta_d)\, p,
\end{equation}
where $K \in \R^{3 \times 3}$ is the \emph{calibration matrix} of the camera, which depends on the its intrinsic parameters, and:
\begin{equation}
  \Pi_0 = \sqmatrixs{1 & 0 & 0 & 0\\ 0 & 1 & 0 & 0\\ 0 & 0 & 1 & 0} \in \R^{3 \times 4}
\end{equation}
is the \emph{standard projection matrix}, which projects 3-D points into the 2-D camera coordinates.
More details about this camera model can be found in Chapter~3 of~\cite{Ma2004}.

\subsection{Closed-Loop Architecture}

As explained in Section~\ref{sec:intro}, we implement a closed-loop system consisting of an observer, which estimates the attacker's intentions using a Receding Horizon Estimation algorithm, and a controller, which plans the path for the robotic arm to defend from the attacker's actions using Receding Horizon Control.
Both, controller and observer, work in a tight coordination.
Even though the controller's main objective is to defend the set $S_{\text{target}}$, its objective also includes a term to force the camera to maintain the sword within its field of view.
Similarly, even though we could have designed an observer that simply estimates the future positions of the attacker's sword, instead we designed an observer that describes the attacker's movements using a fully actuated arm model, forcing our predictions to satisfy the natural constraints of an arm.

\begin{figure}[tp]
  \centering
  \includegraphics[width=.55\linewidth]{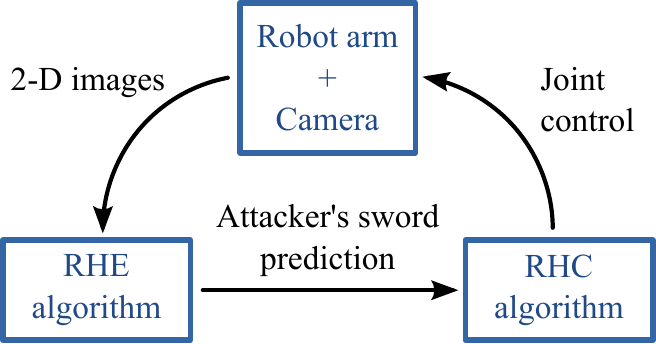}
  \caption{Block diagram of the information flow in the closed-loop system.}
  \label{fig:closed_loop}
\end{figure}

Figure~\ref{fig:closed_loop} shows a block diagram of the closed-loop system.
We explain the details of the observer design in Section~\ref{sec:observer}, and the details of the controller design in Section~\ref{sec:controller}.

\section{Observer Design}
\label{sec:observer}

Once the 2-D camera coordinates of the tip and base of the attacker's sword, $p_t$ and $p_b$, have been obtained via image segmentation, we are faced with the task of extracting useful (or \emph{actionable}) information from this stream of coordinates.
To achieve this goal, our approach is based on the following principles:
first, predictions of the future positions of the sword must account for the non-holonomic constraints of the arm holding it;
second, due to the fast nature of the movements in fencing, it is more useful to predict intentions than exact positions.

We designed an observer based on a Receding Horizon Estimation (RHE) algorithm, which given a time horizon $T_{\text{obs}} > 0$ and a time series on the interval $[t-T_{\text{obs}}, t]$ of 2-D camera coordinates previously collected, estimates the initial angles and the angular velocities of a kinematic arm model we arbitrarily assign to the attacker's arm.
Even though our opponent is assumed to be human, and thus our kinematic model might not match the actual parameters of the attacker's arm, we believe the data produced by it is still better than using a model-less approach to the observation process.
In particular, we chose to use the same 6-joint model we use for the robotic arm to describe the movements of the attacker, noting that if this model is in practice too inaccurate then it is easy to replace it with a different one.

Let $\Theta_a = \sqmatrix{\theta_{a,1}, \dotsc, \theta_{a,6}}^T$ be the joint angles of the attacker's arm, satisfying the following ODE: $\dot{\Theta}_a(t) = u_a(t)$, where $u_a(t) \in \R^6$ is the vector of angular velocities of the joints.
We assume that $G_a(\Theta_a) \in SE(3)$, the transformation taking points from attacker's sword coordinates to world coordinates, is known.
Note that by assuming that $G_a$ is known, we are implicitly assuming that we know the exact location of the attacker's body.
We currently assume that the attacker's body is static.
Also, let $\obs{p}_t(t)$ and $\obs{p}_b(t)$ be the time series of observations in 2-D camera coordinates of the tip of the attacker's sword tip and base, respectively.
Hence, at time $t$ we solve the following optimization problem:
\begin{equation}
  \label{eq:rhe}
  \begin{aligned}
    \min_{\xi, u_a}\,
    &\int_{t-T_{\text{obs}}}^t \hspace{-1em}\lambda(\tau) \!\!\sum_{j \in \set{t,b}}\norm{\obs{p}_j(\tau) - \obs{q}_j(\tau)}_2^2
    + \alpha_1\, \norm{u_a(\tau)}_2^2\, \diff{\tau}\\
    \text{s.t.}\
    &\Theta_a(\tau) \in \sqparen{-\Theta_{\text{max}},\Theta_{\text{max}}},\
    u_a(\tau) \in \sqparen{-u_{\text{max}},u_{\text{max}}},\\
    &\dot{\Theta}_a(\tau) = u_a(\tau),\
    \Theta_a(t-T_{\text{obs}}) = \xi,\\
    &q_j(\tau) = G_{\text{cam}}^{-1}\parenb{\Theta_d(\tau)}\, G_a\parenb{\Theta_a(\tau)}\, p_j,\\
    &\sqparen{q_j(\tau)}_{\axis{z}_{\text{cam}}}\, \obs{q}_j(\tau) = \, K\, \Pi_0\, q_j(\tau),\\
    &\sqparen{q_j(\tau)}_{\axis{z}_{\text{cam}}} \geq \epsilon_{\text{min}},\
    \forall j \in \set{t,b},\
    \forall \tau \in [t-T_{\text{obs}},t],\\
  \end{aligned}
\end{equation}
where $\lambda(\tau)$ is a forgetfulness factor putting more weight on recent observations, defined by $\lambda(\tau) = \exp\parenb{\alpha_2\, (\tau-t)}$, $p_t$ and $p_b$ are defined as in equation~\eqref{eq:pt_pb}, $\alpha_1, \alpha_2 > 0$ are constant parameters, and $\epsilon_{\text{min}}$ is the closest distance an object can be from the focal plane of the camera.

We solve the optimal estimation problem in equation~\eqref{eq:rhe}, at regular intervals of $T_{\text{int}} > 0$ seconds, by transforming this problem into a nonlinear programming problem via time discretization as described in Chapter~4 of~\cite{Polak1997}.
As part of the discretization process, we assume that $u_a(t)$ is constant over the interval $[t-T_{\text{obs}},t]$.
The intuition behind this simplification comes from the fencing rules~\cite{fencingrules}, which require an attack to be a continuous motion towards the target.
This assumption is also consistent with the speed of fencing movements, which are completed in the order of hundreds of milliseconds; in an interval of that length we can only acquire and process a few tens of images, hence it is better to use a simple movement model than running the risk of over-fitting the data.

\section{Controller Design}
\label{sec:controller}

In fencing, an attack can be interrupted by contact between the two fencers' swords~\cite{fencingrules}. The most basic form of blocking consists of the defender placing his sword between the attacker's sword and the target area.
Since the attacker cannot interrupt his motion towards the target area without ending his attack, this is often sufficient to achieve a successful block.
We designed a Receding Horizon Control (RHC) algorithm stop the attacker's sword with the defender's before the attacker reaches any of the points in $S_{\text{target}}$.
The RHC algorithm works by iteratively solving an optimal control problem, with time horizon $T_{\text{con}} > 0$, every $T_{\text{int}} > 0$ seconds, and it achieves the desired objective after a careful design of the objective and constraint functions, as explained below.

We define the \emph{blocking plane} as the plane parallel to the defender's sword with normal vector $\axis{x}_d$.
Thus, given a set joint angles of the defender $\Theta_d \in \R^6$, a point $p \in \R^3$, and the function $h_{\text{bp}}(p,\Theta_d) = \axis{x}_d^T(\Theta_d)\, \parenb{p - G_d(\Theta_d)\, p_b}$, the blocking plane can be described as the set $\set{p \in \R^3 \mid h_{\text{bp}}(p,\Theta_d) = 0}$.
Our controller must maintain the attacker's sword and the set $S_{\text{target}}$ divided by the blocking plane to ensure the safety of the defender.
Mathematically we write the condition above as the following constraints:
\begin{equation}
  \label{eq:hbp_cons}
  \begin{aligned}
    &h_{\text{bp}}\parenb{G_a(\Theta_a)\, p_j, \Theta_d} \geq 0,\ \forall j \in \set{t,b}, \quad \text{and},\\
    &h_{\text{bp}}\parenb{s, \Theta_d} \leq 0,\ \forall s \in S_{\text{target}}.
  \end{aligned}
\end{equation}

\begin{figure}[tp]
  \centering
  \resizebox{!}{.36\linewidth}{%
    \fbox{%
      \clipbox{50pt 15pt 0pt 35pt}{%
        \input{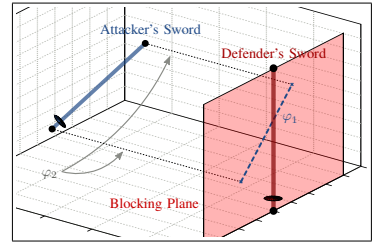}%
      }%
    }%
  }
  \caption{Diagram of the blocking plane (in $\paren{\axis{x}_d,\axis{y}_d,\axis{z}_d}$ coordinates), the projection of the attacker's sword into the blocking plane (with length $\phi_1$, defined in equation~\eqref{eq:phi1}), and the distance of the attacker's sword ends to the blocking plane (their sum is $\phi_2$, defined in equation~\eqref{eq:phi2}).}
  \label{fig:constraints}
\end{figure}

Since the length of the sword is finite, it is not enough to enforce that the attacker's sword is separated from $S_{\text{target}}$. We must also be sure that if the attacker's sword is close to the blocking plane then the defender's sword must be within reach of the attacker's sword.
We achieve this goal by imposing minimum and maximum distances between the defender's sword and the projection of the attacker's sword onto the blocking plane.
Recall that the blocking plane is normal to $\axis{x}_d$.
Thus, given $p \in \R^3$ and $\Theta_d, \Theta_a \in \R^6$, if we define:
\begin{equation}
  \label{eq:hbox_cons}
  h_{\text{box}}(p,\Theta_d,\Theta_a) =
  \begin{pmatrix}
    \sqparen{G_d^{-1}(\Theta_d)\, G_a(\theta_a)\, p - p}_{\axis{y}_d}\\
    \sqparen{G_d^{-1}(\Theta_d)\, G_a(\theta_a)\, p - p}_{\axis{z}_d}\\
  \end{pmatrix},
\end{equation}
then the constraint described above can be written as $b_{\text{min}} \leq h_{\text{box}}(p_j,\Theta_d,\Theta_a) \leq b_{\text{max}}$ for each $j \in \set{t,b}$, where $b_{\text{min}}, b_{\text{max}} \in \R^2$ are the minimum and maximum distances allowed, respectively.
The blocking plane and box constraints are shown in Figure~\ref{fig:constraints}.

To ensure that the defender's sword stops the attacker's when it approaches the blocking plane, we penalize the swords being parallel.
Mathematically, we achieve this goal by including the following term in the objective function of our RHC optimization problem:
\begin{equation}
  \label{eq:phi1}
  \phi_1(\Theta_d,\Theta_a) = \abs{\axis{z}_d^T(\Theta_d)\, \parenb{G_a(\Theta_a)\, p_t - G_a(\Theta_a)\, p_b}},
\end{equation}
which is minimized when both swords are orthogonal.
Note that even though $\phi_1$ is a non-smooth function, it can be easily implemented as part of a nonlinear programming problem using the traditional epigraph transformation (e.g., see Chapter~4.1.3 in~\cite{Boyd2004}).

We also penalize the possibility that the swords could avoid one another if the projection of the attacker's sword onto the blocking plane is very small, i.e., if it is approaching in a ``stabbing'' motion.
Consider the following function calculating the sum of the distances between the ends of the attacker's sword and the blocking plane:
\begin{equation}
  \label{eq:phi2}
  \phi_2(\Theta_d,\Theta_a) = \sum_{j \in \set{t,b}}
  \abs{\axis{x}_{d}^T(\Theta_d)\, \parenb{G_a(\Theta_a)\, p_j - G_d(\Theta_d)\, p_b}}.
\end{equation}
Minimizing this function reduces the distances between both ends of the opponent's sword and the blocking plane.
Note that when $\phi_2$ is minimized, the attacker's sword becomes parallel to the blocking plane, and at the same time the projection of the sword onto the blocking plane is maximized, thus removing the possibility of a ``stabbing'' motion.
This function also helps ensure that the swords will eventually come into contact, since aim to reduce the overall distance between both swords.
Also note that thanks to the constraint in equation~\eqref{eq:hbp_cons} we can simply remove the absolute values from the formula of $\phi_2$.

Note that all the functions above involve the RHC algorithm using information provided by the observer to estimate the angles of the attacker's joints, $\Theta_a$.
However, the quality of the estimations provided by the observer are directly related to the data captured by the camera, which in turn is related to the angles of the defender's joints, $\Theta_d$.
Hence, we add a term to the objective function of the RHC optimization problem that is minimized when the attacker's sword is centered in the 2-D camera coordinate frame:
\begin{equation}
  \label{eq:phi3}
  \phi_3(\Theta_d, \Theta_a) =
  \norm{%
    \frac{K\, \Pi_0\, G_{\text{cam}}^{-1}(\Theta_d)\, G_a(\Theta_a)\, p_b}%
    {\sqmatrix{G_{\text{cam}}^{-1}(\Theta_d)\, G_a(\Theta_a)\, p_b}_{\axis{z}_{\text{cam}}}}%
  }_2^2.
\end{equation}
Although it would be ideal to implement a \emph{hard} constraint requiring the controller to keep the opponent's sword in the camera's view at all times, in practice this leads to an over-constrained problem which often has no feasible solutions.
Therefore, we implement a \emph{soft} constraint adding $\phi_3$ to the running cost.

Recall that at time $t$ we know the defender's joint angles $\Theta_d(t)$, an estimation of the attacker's joint angles $\Theta_a(t)$, and an estimation of the attacker's joint angular velocities $u_a(t)$, the first from angle encoders mounted on the robotic arm, and the last two as results of the observer presented in Section~\ref{sec:observer}.
For prediction purposes, we assume that the attacker will apply constant angular velocities to all its joints on the interval $\sqparen{t,t+T_{\text{con}}}$ equal to $u_a(t)$.

Combining the functions defined above, the results from the observer, together with dynamics constraints and a regularization term in the objective, yields the following optimal control problem:
\begin{equation}
  \label{eq:rhc}
  \begin{aligned}
    \min_{u_d}\,
    &\int_{t}^{t_f} \!\!\paren{\gamma_3\, \phi_3\parenb{\Theta_d(\tau), \Theta_a(\tau)}
    + \gamma_4\, \norm{u_d(\tau)}_2^2} \, \diff{\tau} +\\
    &\hspace{10pt}+ \gamma_1\, \phi_1\parenb{\Theta_d(t_f),\Theta_a(t_f)}
    + \gamma_2\, \phi_2\parenb{\Theta_d(t_f),\Theta_a(t_f)},\\
    \text{s.t.}\
    &\Theta_d(\tau) \in \sqparen{-\Theta_{\text{max}},\Theta_{\text{max}}},\
    u_d(\tau) \in \sqparen{-u_{\text{max}},u_{\text{max}}},\\
    &\dot{\Theta}_d(\tau) = u_d(\tau),\
    \dot{\Theta}_a(\tau) = u_a(t),\\
    &h_{\text{bp}}\parenb{G_a\parenb{\Theta_a(\tau)}\, p_j, \Theta_d(\tau)} \geq 0,\\
    &h_{\text{bp}}\parenb{s, \Theta_d(\tau)} \leq 0,\\
    &b_{\text{min}} \leq h_{\text{box}}\parenb{p_j, \Theta_d(\tau), \Theta_a(\tau)} \leq b_{\text{max}},\\
    &\forall s \in S_{\text{target}},\ \forall j \in \set{t,b},\ \forall \tau \in [t,t_f],
  \end{aligned}
\end{equation}
where $t_f = t + T_{\text{con}}$, and $\gamma_1, \gamma_2, \gamma_3, \gamma_4 > 0$ are scaling parameters.
Similar to the optimization problem in equation~\eqref{eq:rhe}, we solve this optimization problem every $T_{\text{int}}$ seconds by discretizing it using the techniques described in Chapter~4 in~\cite{Polak1997}.

\section{Simulations}
\label{sec:examples}

We simulated the behavior of our closed-loop scheme for a collection of predetermined movements of the attacker.
We modeled both defender and attacker arms based on a FANUC robotic arm model LRMate~200iC~\cite{fanuc_lrmatic_200ic}.
Also, we modeled the camera mounted on the defender's arm on a Point~Grey camera model Flea~3~\cite{ptgrey_flea3}.
We discretized both RHE and RHC optimization problems, described in equations~\eqref{eq:rhe} and~\eqref{eq:rhc}, with a step-size $T_{\text{disc}}$, and we solved the resulting discretized nonlinear programming optimization problem using the SNOPT library~\cite{Gill2008}.
The simulations were programmed using the language Python, and they were executed in a two-processor Xeon~E5-2680 computer running at $2.7\, \sqparen{\textrm{GHz}}$.
A detailed list of all the parameters used in our simulations is shown in Table~\ref{tab:params}.
Note that the camera calibration matrix is $K = \textrm{diag}\paren{f_\ell\, \sigma_p, f_\ell\, \sigma_p, 1}$, where $f_\ell$ is the camera focal length, and $\sigma_p$ is the camera pixel size.
Also note that even though the sword length is $880\, \sqparen{\textrm{mm}}$, the upper limit for the box constraint in equation~\eqref{eq:hbox_cons} along the $\axis{z}_d$-axis is $600$, since this value reflects the fact that attacks are hard to stop in practice with the tip of the defender's sword.

\begin{table}[tp]
  \caption{%
    List of the parameters used in the simulations.
  }
  \label{tab:params}
  \centering
  \begin{tabular}{c|c}
    Parameter $\sqparen{unit}$ & Value \\
    \hline
    Robot segment lengths $\sqparen{\textrm{mm}}$ & $\sqmatrixs{50 & 250 & 300 & 75 & 80}$ \\
    Sword length $\sqparen{\textrm{mm}}$ & 880 \\
    $\Theta_{\text{max}}\, \sqparen{^\circ}$ & $\sqmatrixs{180 & 150 & 190 & 190 & 220 & 360}^T$ \\
    $\Theta_{0}\, \sqparen{^\circ}$ & $\sqmatrixs{0 & 135 & 95 & 0 & 115 & -15}^T$ \\
    $u_{\text{max}}\, \sqparen{^\circ/\textrm{s}}$ & $\sqmatrixs{350 & 350 & 400 & 450 & 450 & 720}^T$ \\
    Camera focal length $f_\ell\, \sqparen{\textrm{mm}}$ & $12.5$ \\
    Camera pixel size $\sigma_p\, \sqparen{\mu\textrm{m}}$ & $4.5$ \\
    Camera resolution $\sqparen{\textrm{pixels}}$ & $1600 \times 1200$ \\
    $T_{\text{con}}\, \sqparen{\textrm{ms}}$ & $100$ \\
    $T_{\text{obs}}\, \sqparen{\textrm{ms}}$ & $50$ \\
    $T_{\text{int}}\, \sqparen{\textrm{ms}}$ & $15$ \\
    $T_{\text{disc}}\, \sqparen{\textrm{ms}}$ & $5$ \\
    $\alpha_1$ & $10^{-2}$ \\
    $\alpha_2$ & $10^{-2}$ \\
    $\gamma_1$ & $10^{-2}$ \\
    $\gamma_2$ & $1$ \\
    $\gamma_3$ & $10^{-4}$ \\
    $\gamma_4$ & $10^{-4}$ \\
    $b_{\text{min}}$ & $\sqmatrixs{-500 & 10}^T$ \\
    $b_{\text{max}}$ & $\sqmatrixs{500 & 600}^T$ \\
  \end{tabular}
\end{table}

We include simulations under two different scenarios below.
In the first case the attacker's joint angular velocities, $u_a(t)$, are constant over the whole , the second where $u_a(t)$ changes over time.
We also present results validating the use of the penalty function $\phi_3$, designed to keep the attacker's sword within the field of view of the defender's camera, in the RHC optimization problem in equation~\eqref{eq:rhc}.

\subsection{Example 1: $u_a$ is constant}
\label{subsec:example1}

We begin with a sample trajectory where the assumption that $u_a(t)$ is constant holds.
The initial joint positions for the robot and the opponent are given by $\Theta_d(0) = \Theta_a(0) = \Theta_0$, with $\Theta_0$ given in Table~\ref{tab:params}, and the opponent's trajectory is defined by:
\begin{equation}
  u_a(t) = \sqmatrixs{-21 & -210 & -400 & 345 & -63 & 0}^T.
\end{equation}
This trajectory defines an attack to the robot's left side.

\begin{figure*}
  \centering
  \begin{tabular}{ccccc}
    \topfig{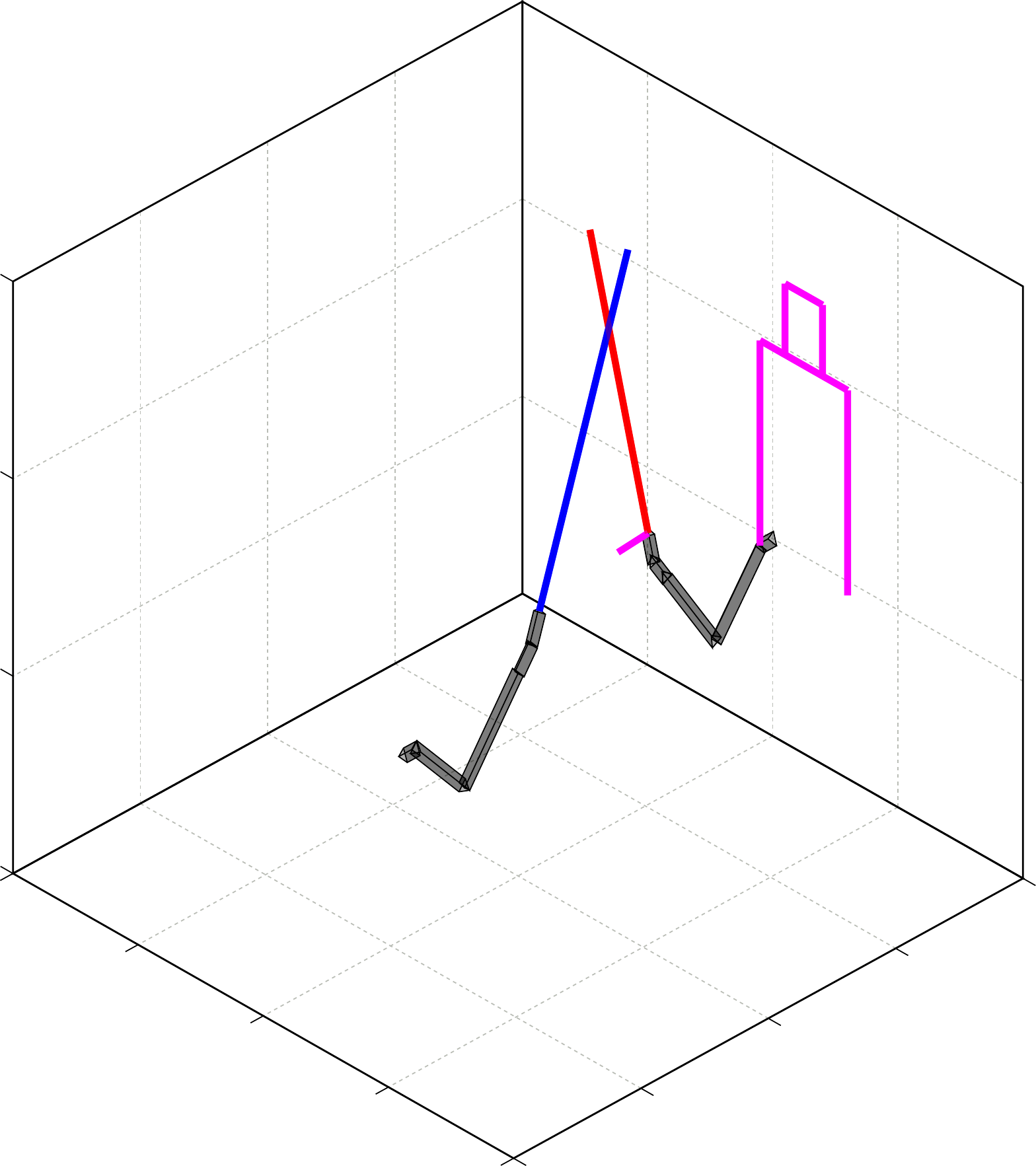} &
    \topfig{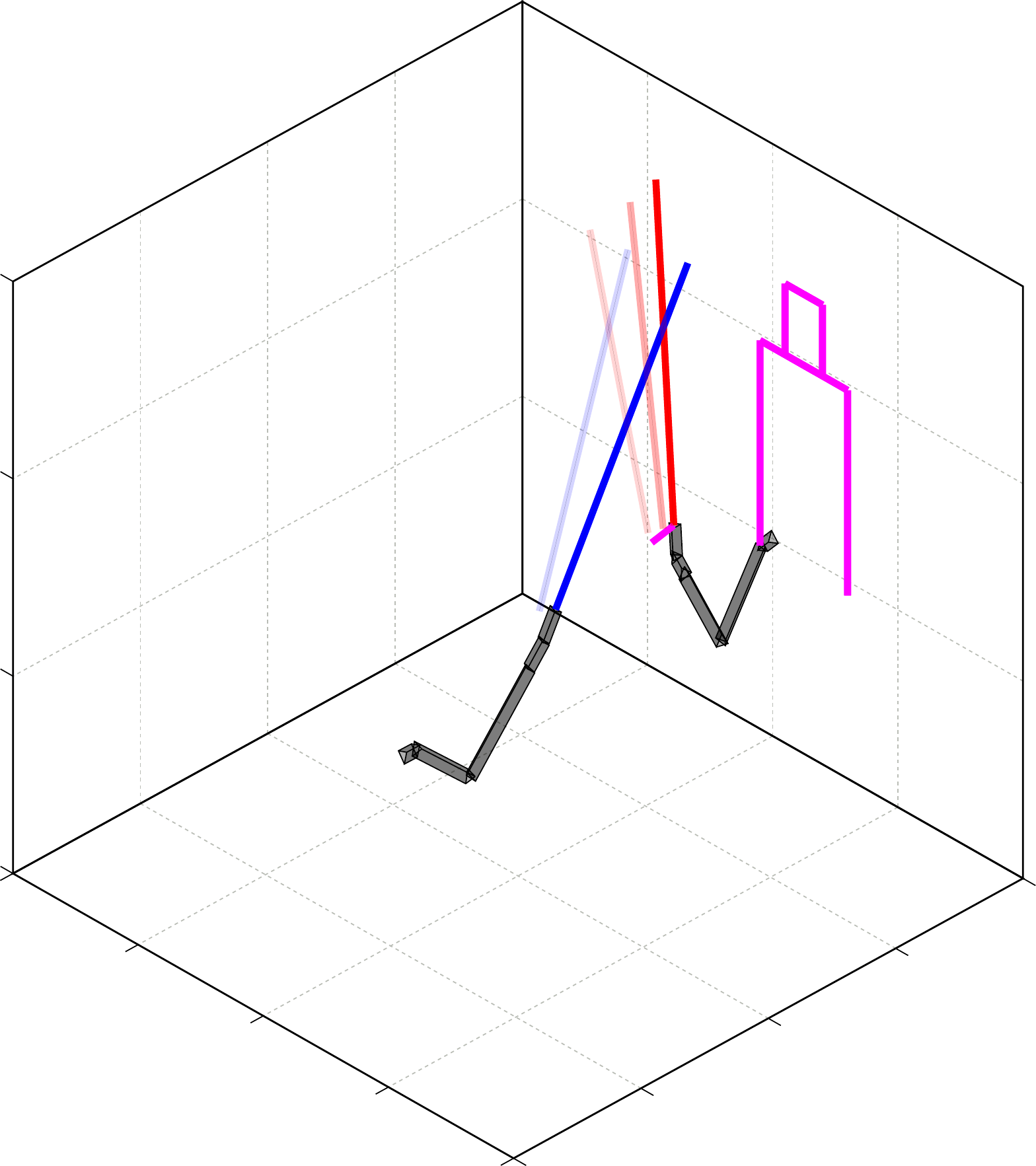} &
    \topfig{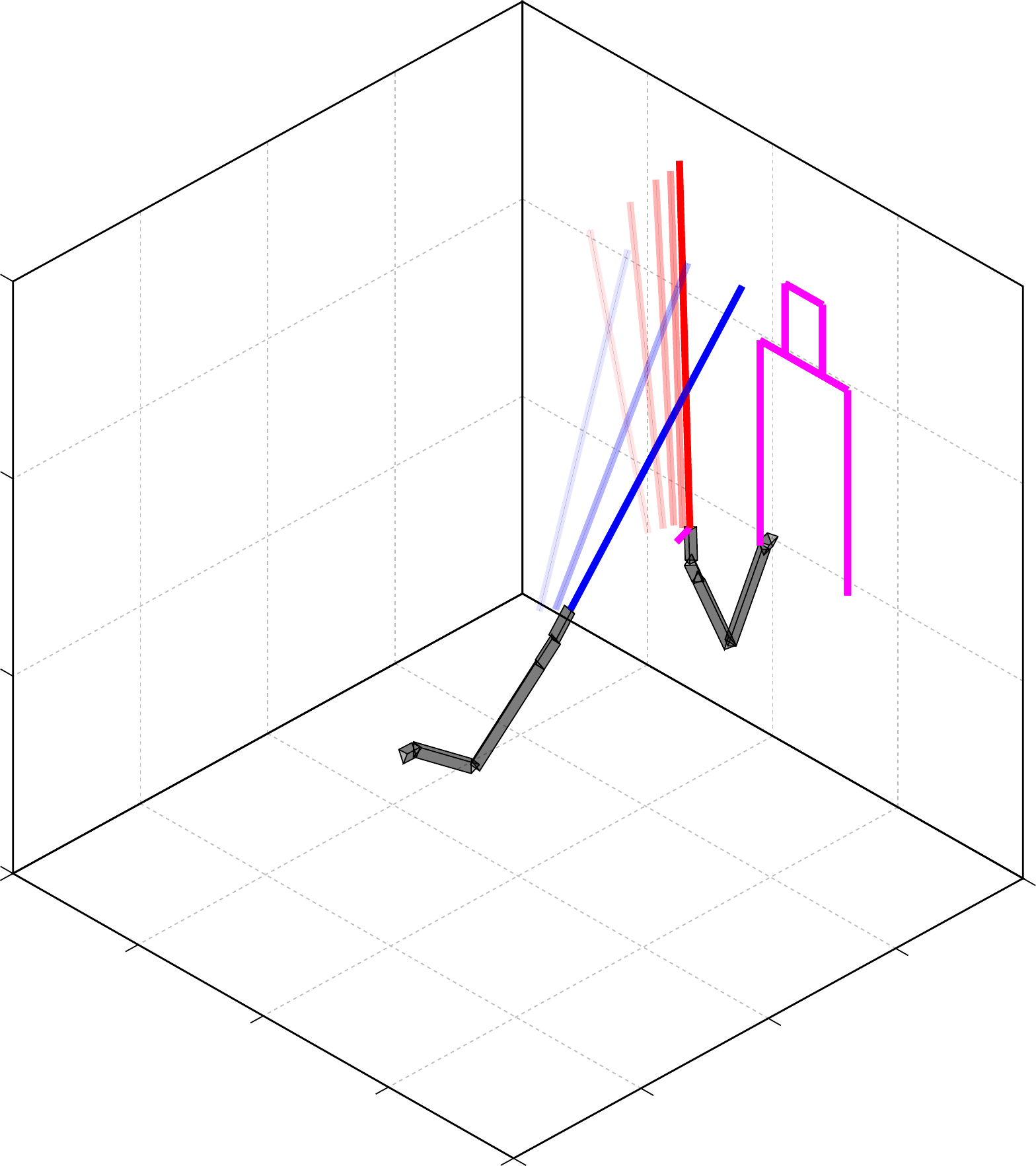} &
    \topfig{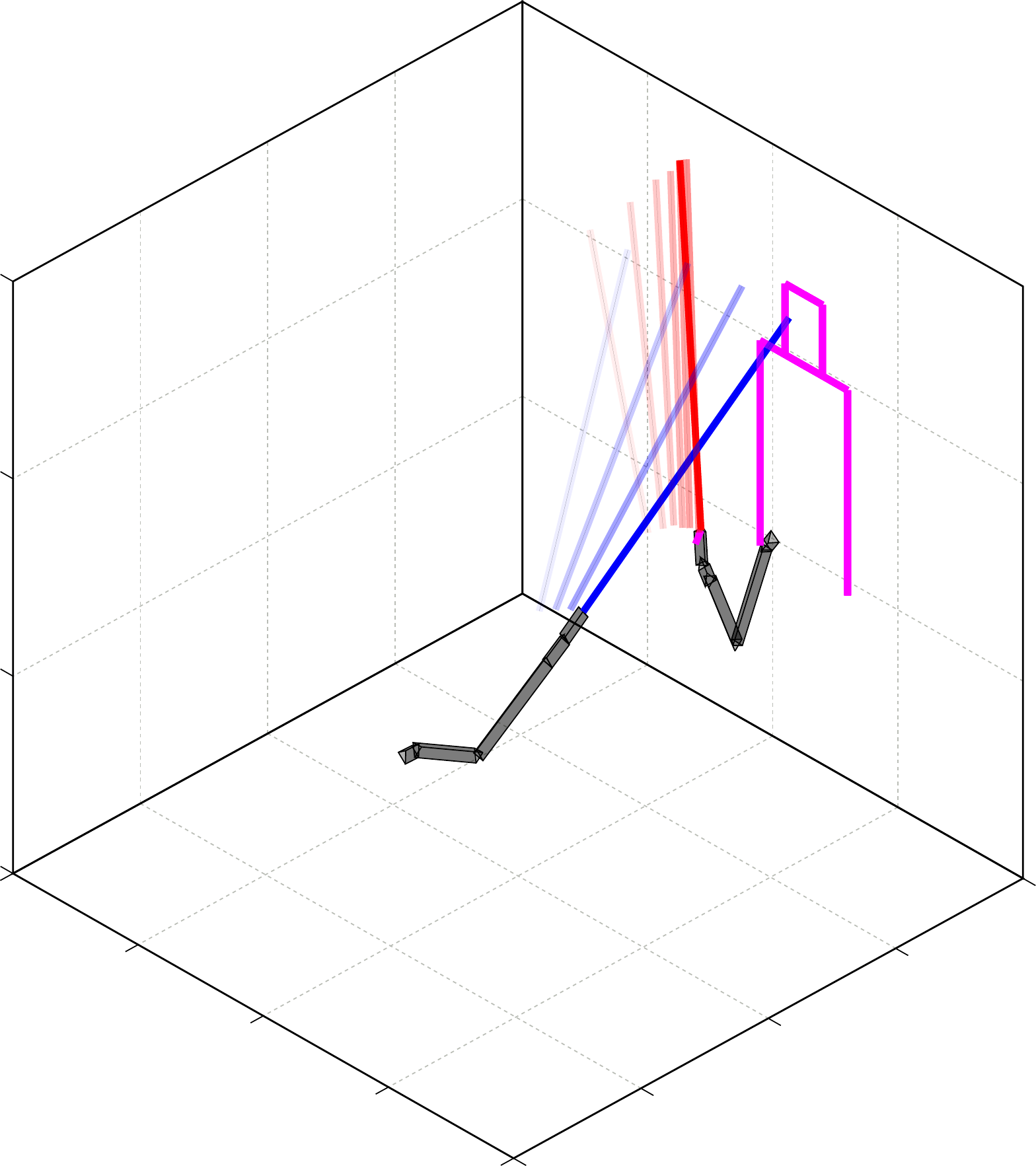} &
    \topfig{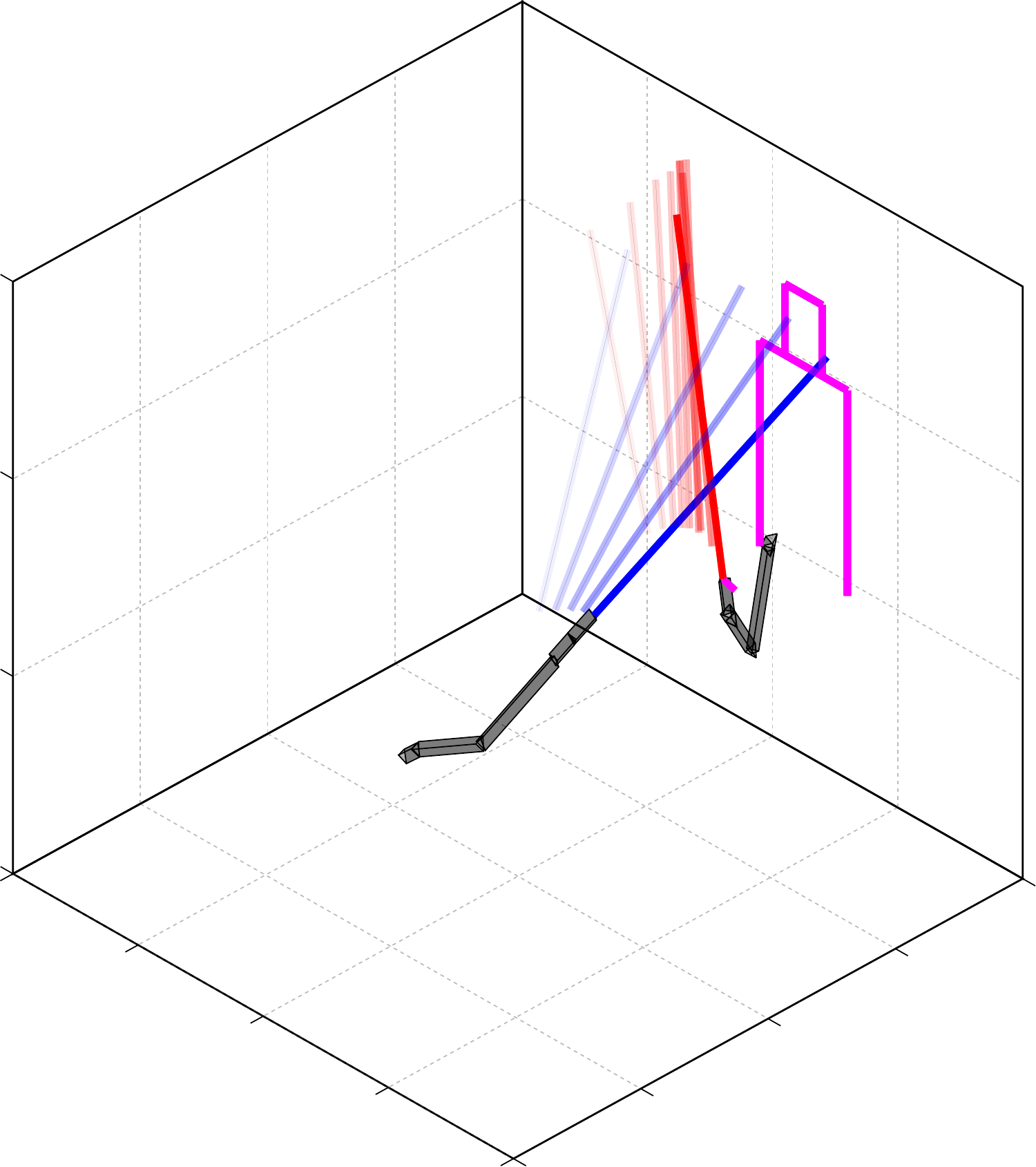} \\[5pt]
    \botfig{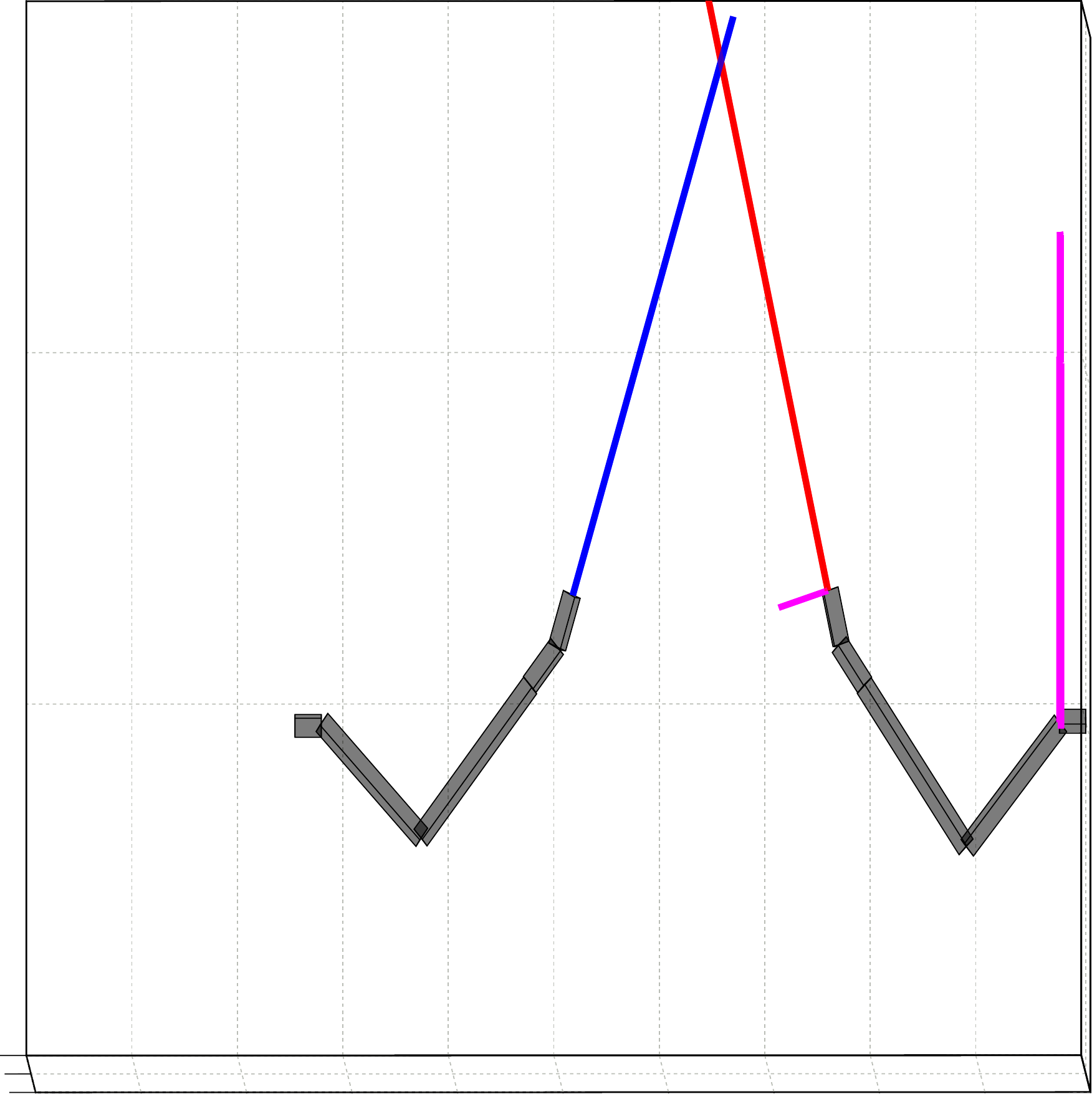} &
    \botfig{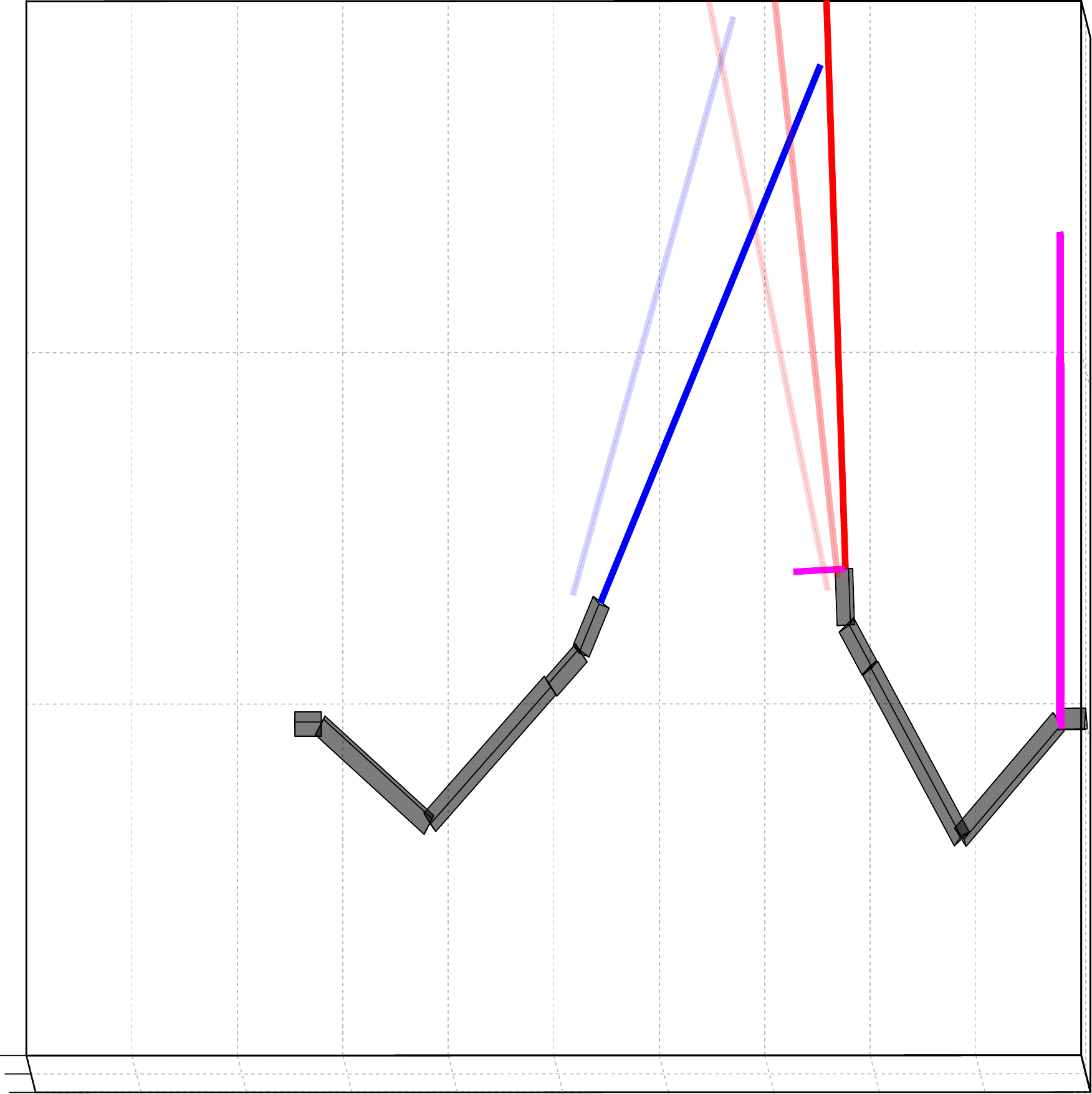} &
    \botfig{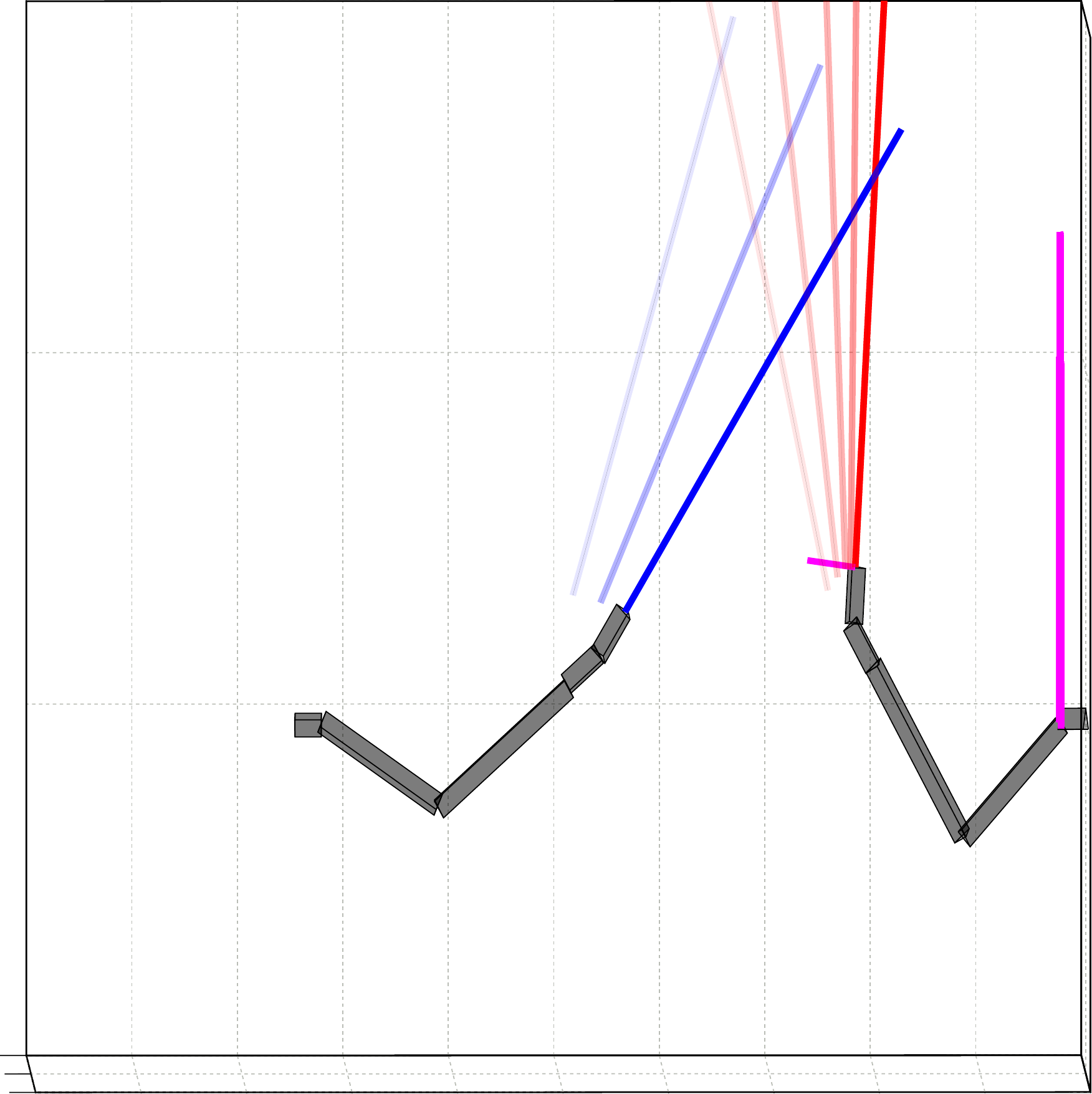} &
    \botfig{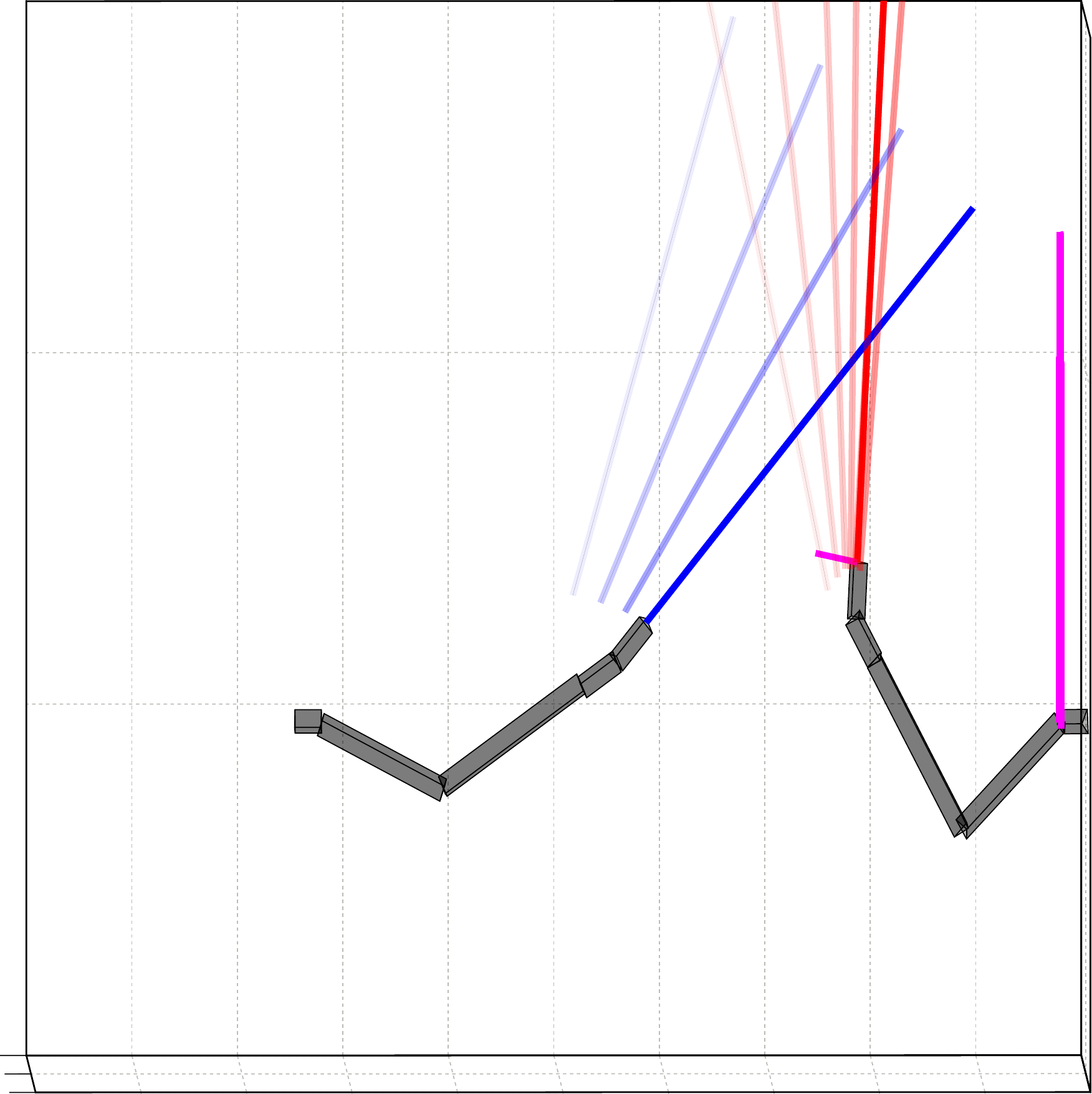} &
    \botfig{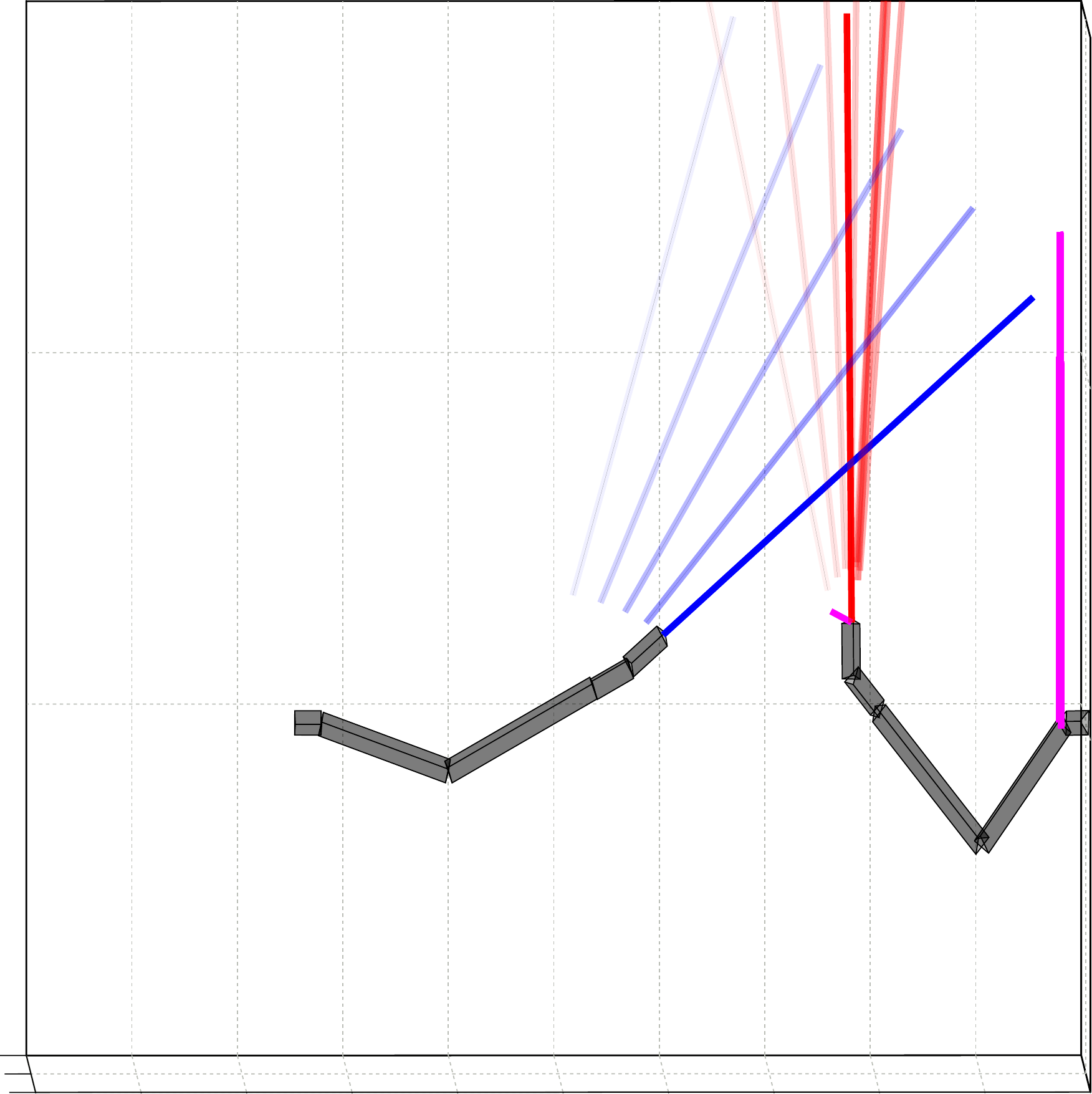} \\[3pt]
    $t = 0\,\sqparen{\textrm{ms}}$ &
    $t = 36\,\sqparen{\textrm{ms}}$ &
    $t = 72\,\sqparen{\textrm{ms}}$ &
    $t = 108\,\sqparen{\textrm{ms}}$ &
    $t = 145\,\sqparen{\textrm{ms}}$ \\
  \end{tabular}
  \caption{%
    Snapshots of the simulation in Section~\ref{subsec:example1}.
    The attacker's sword is shown in blue, the defender's sword in red, and the set $S_\text{target}$ is shown in purple.
    Each column corresponds to a different time, the top row from an inclined point of view, the bottom row from a lateral point of view.
  }
  \label{fig:ex1_trajectory}
\end{figure*}

Figure~\ref{fig:ex1_trajectory} shows the optimal trajectory computed by our closed-loop algorithms for this scenario.
The robot successfully blocks the opponent's attack by achieving blade contact at time $t_f = 145\, \sqparen{\textrm{ms}}$.


\subsection{Example 2: $u_a$ is not constant}
\label{subsec:example2}

\begin{figure*}
  \centering
  \begin{tabular}{ccccc}
    \topfig{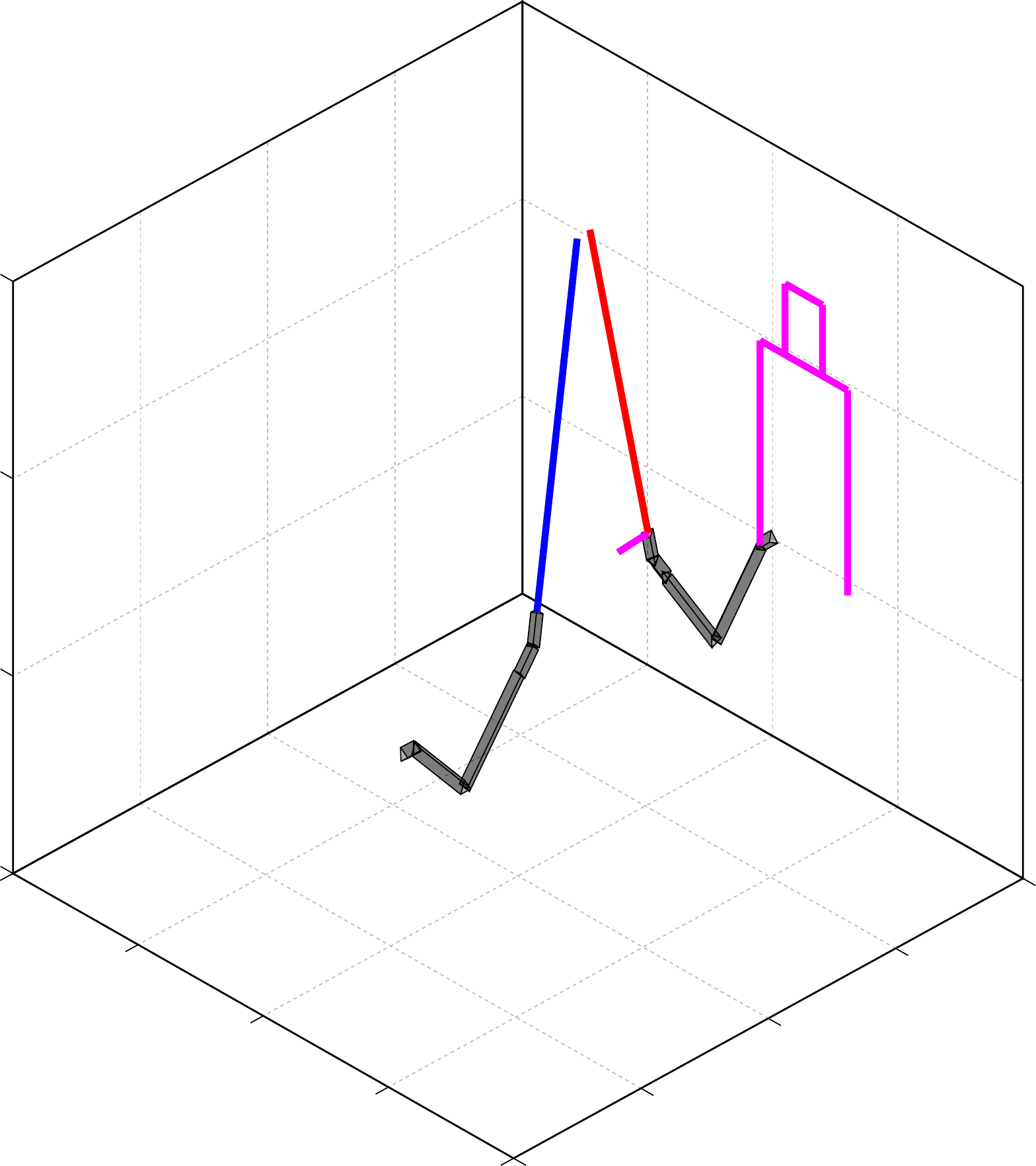} &
    \topfig{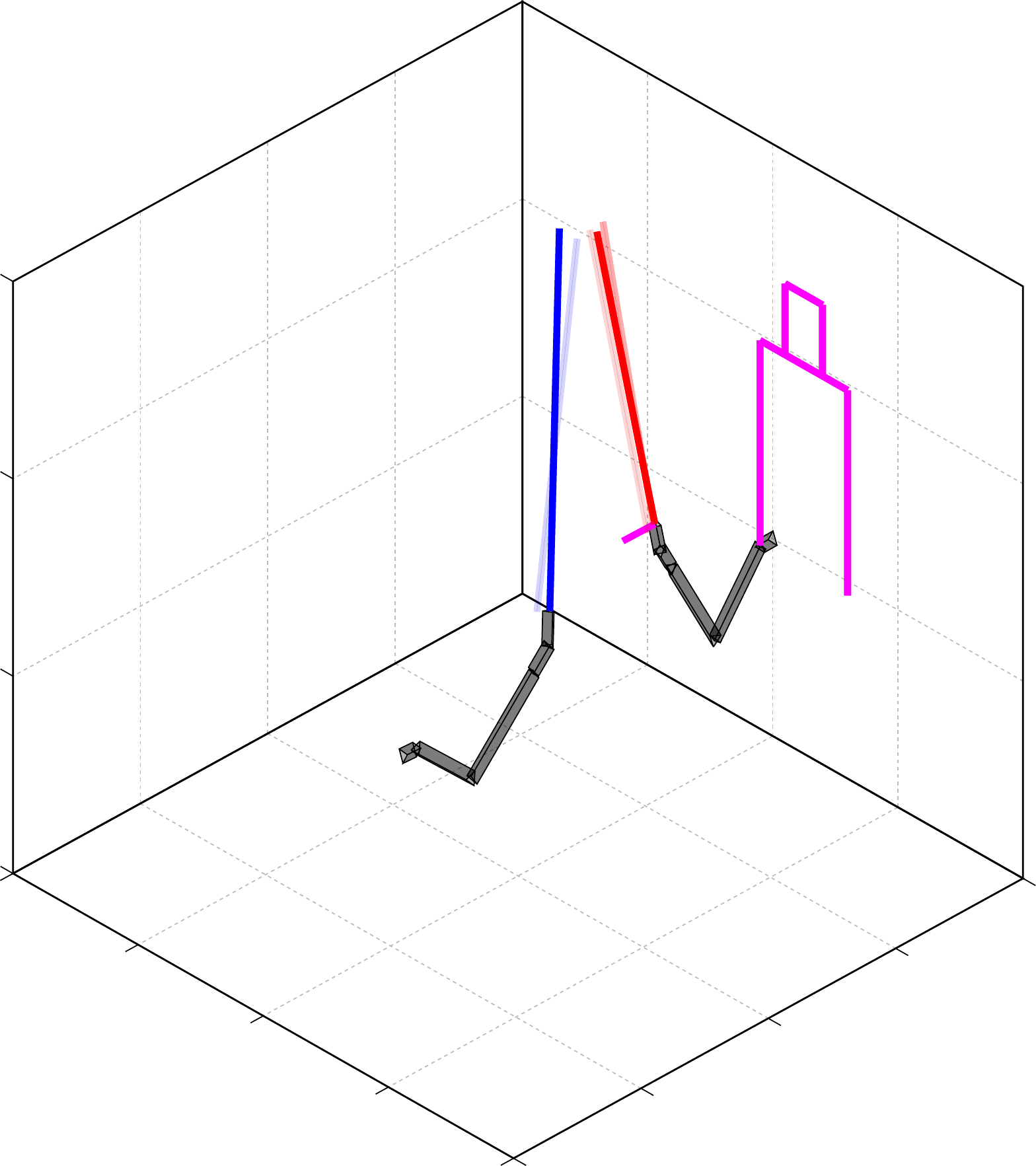} &
    \topfig{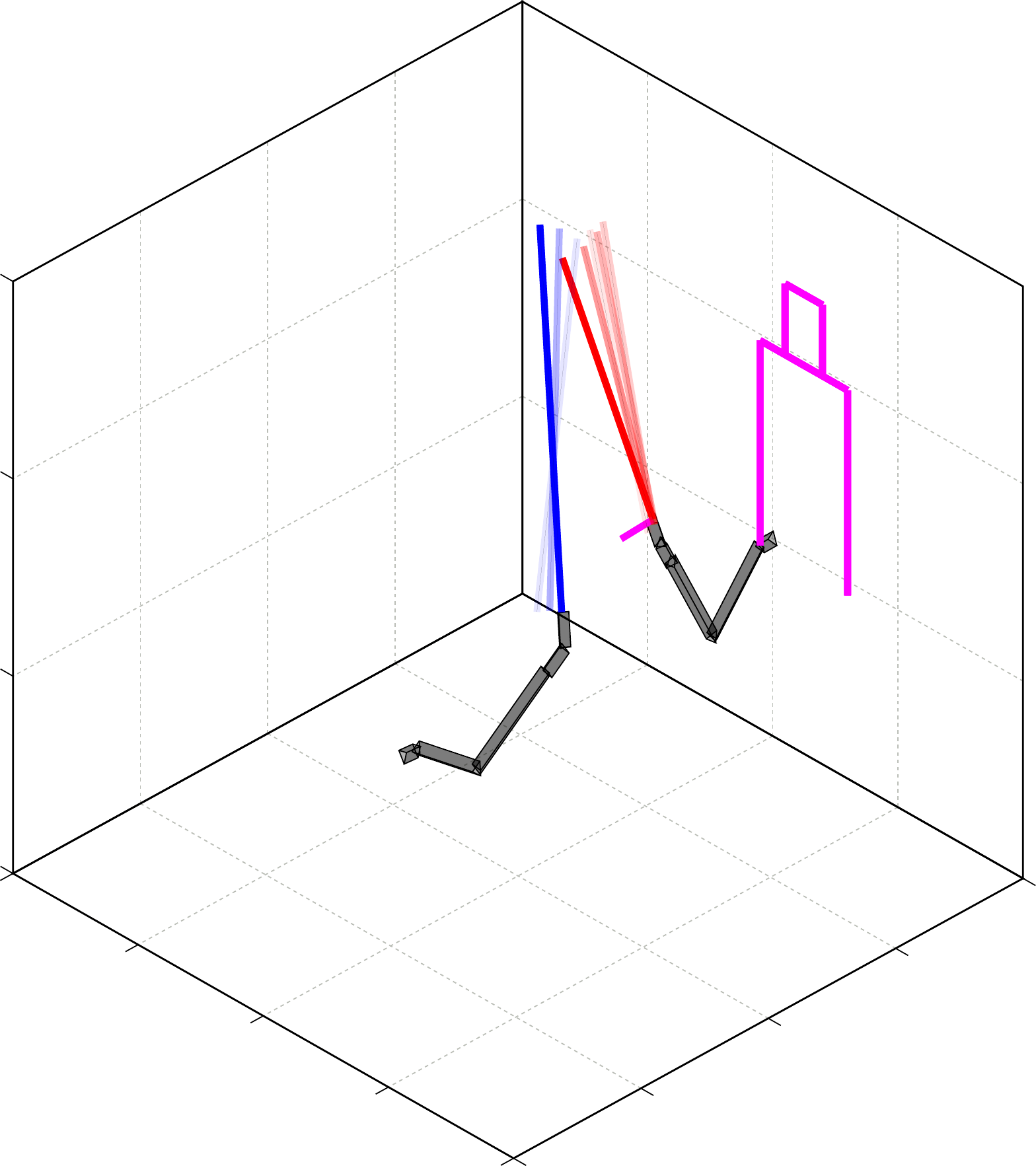} &
    \topfig{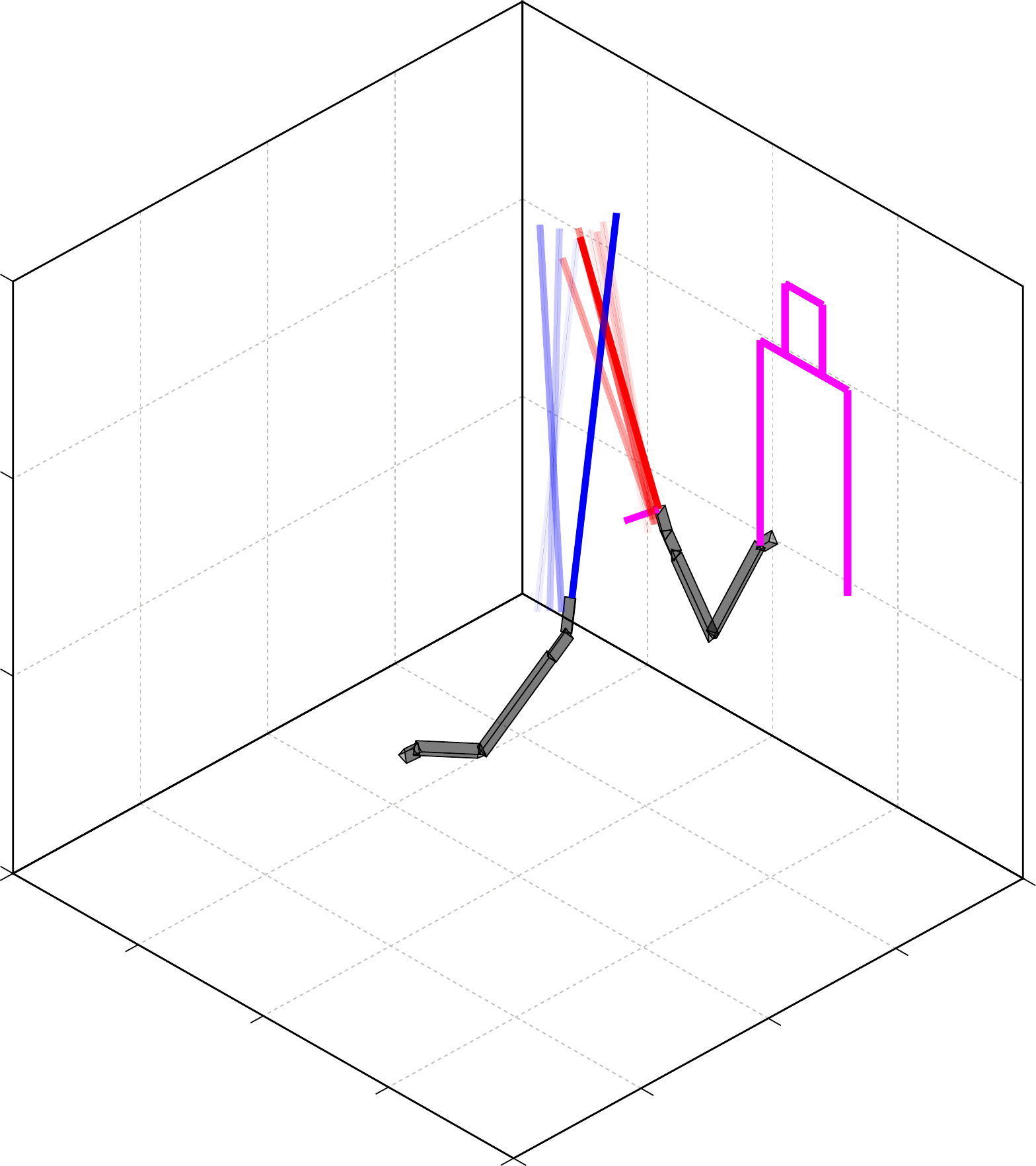} &
    \topfig{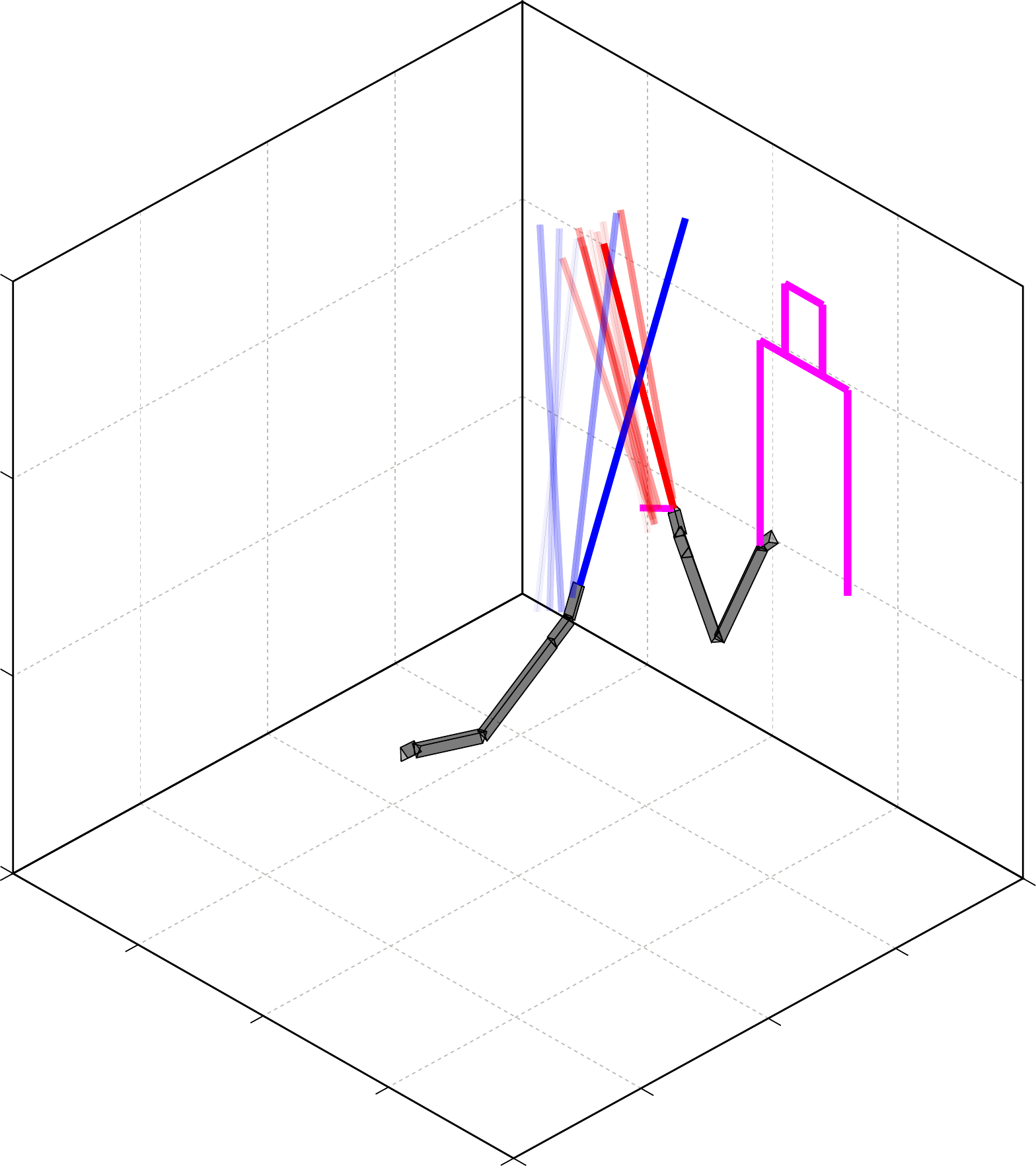} \\[5pt]
    \botfig{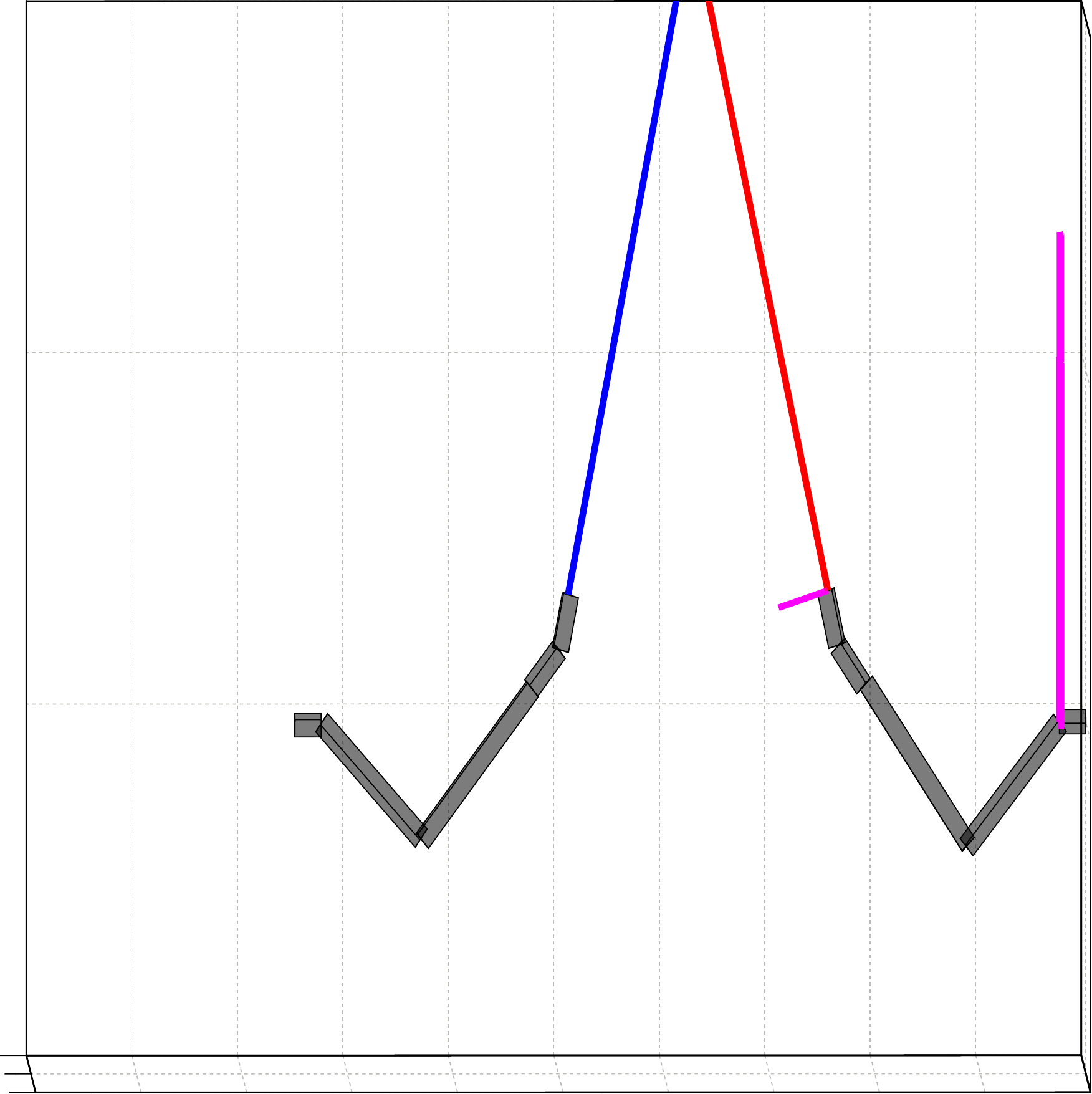} &
    \botfig{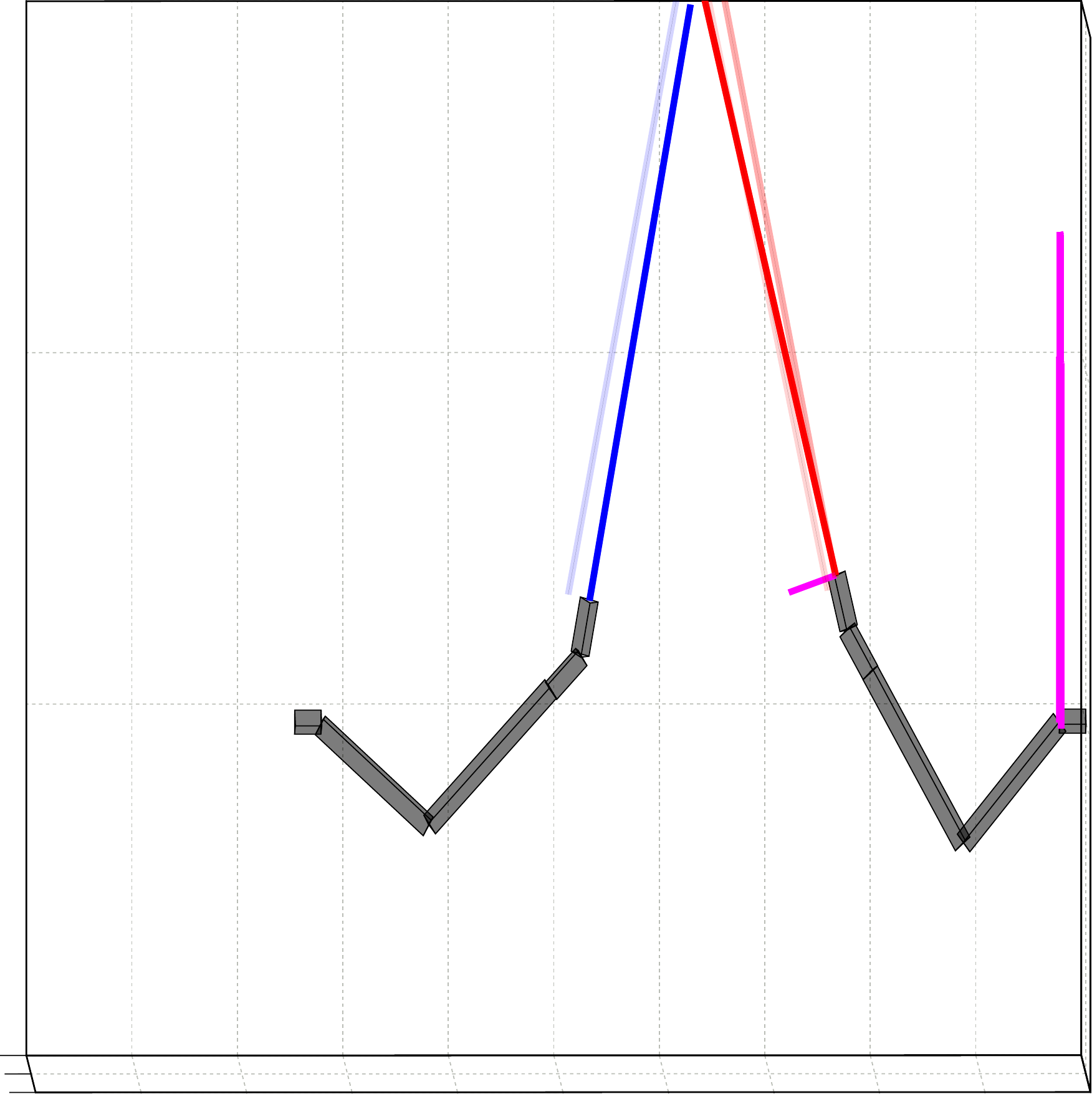} &
    \botfig{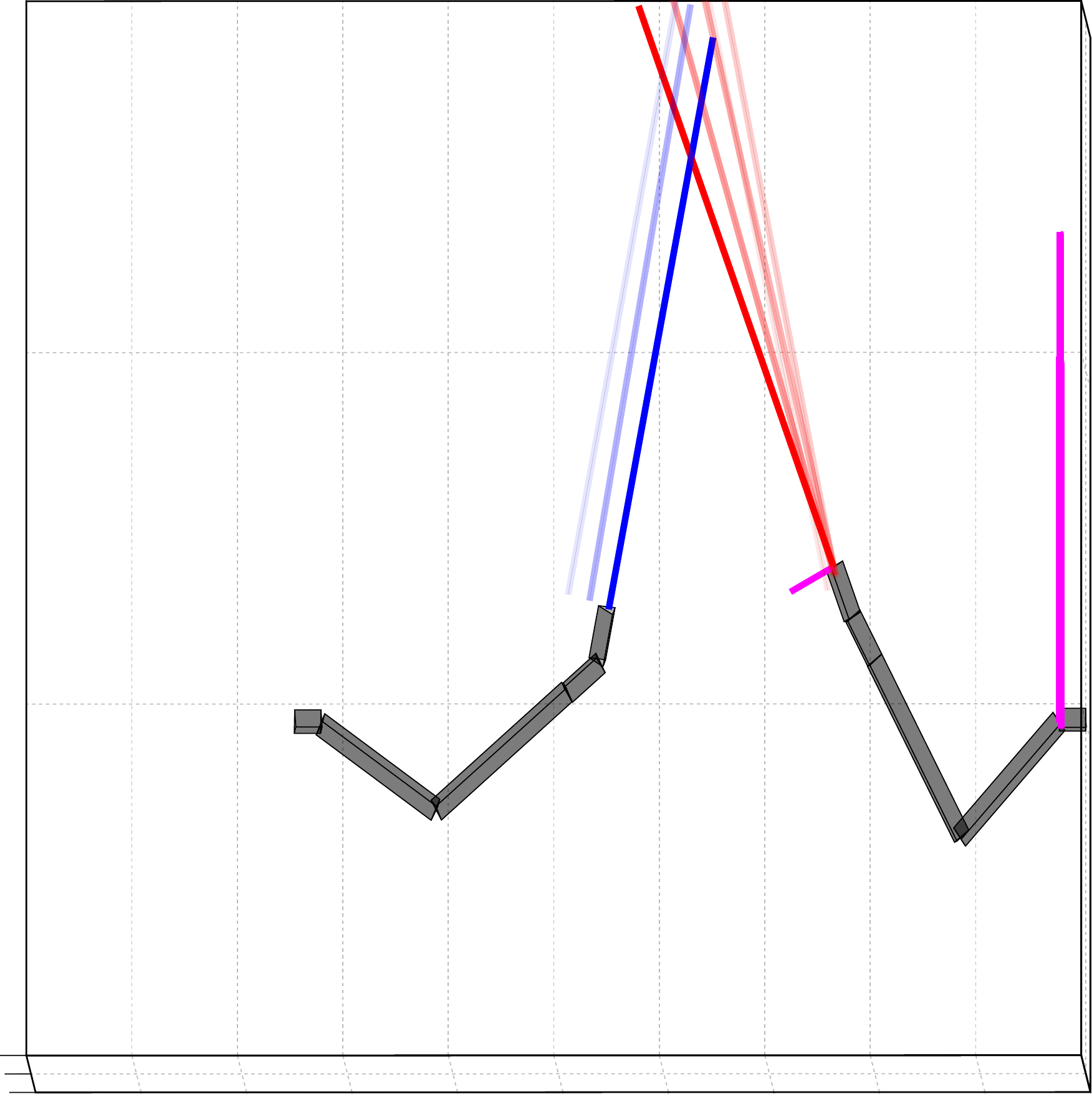} &
    \botfig{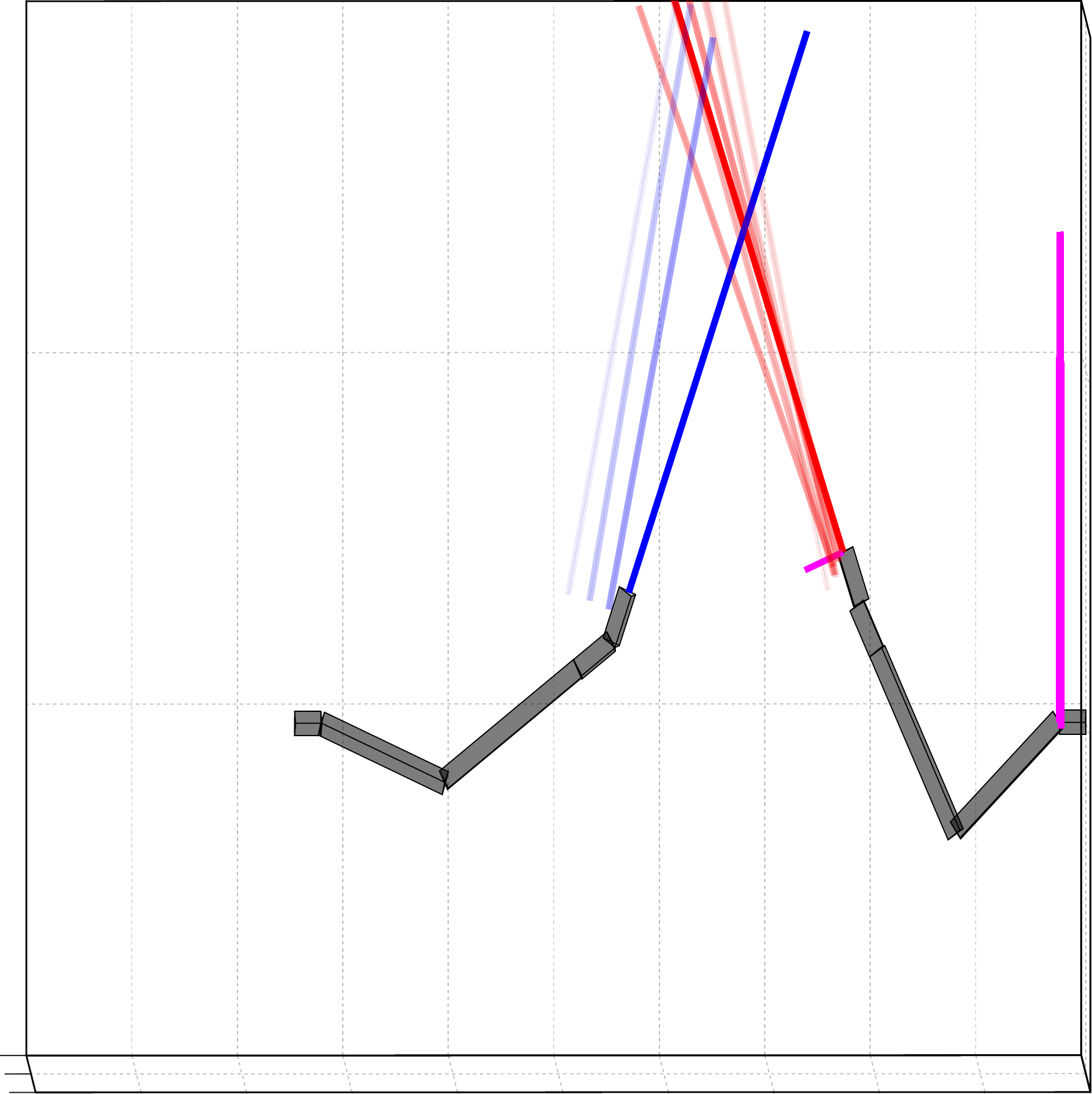} &
    \botfig{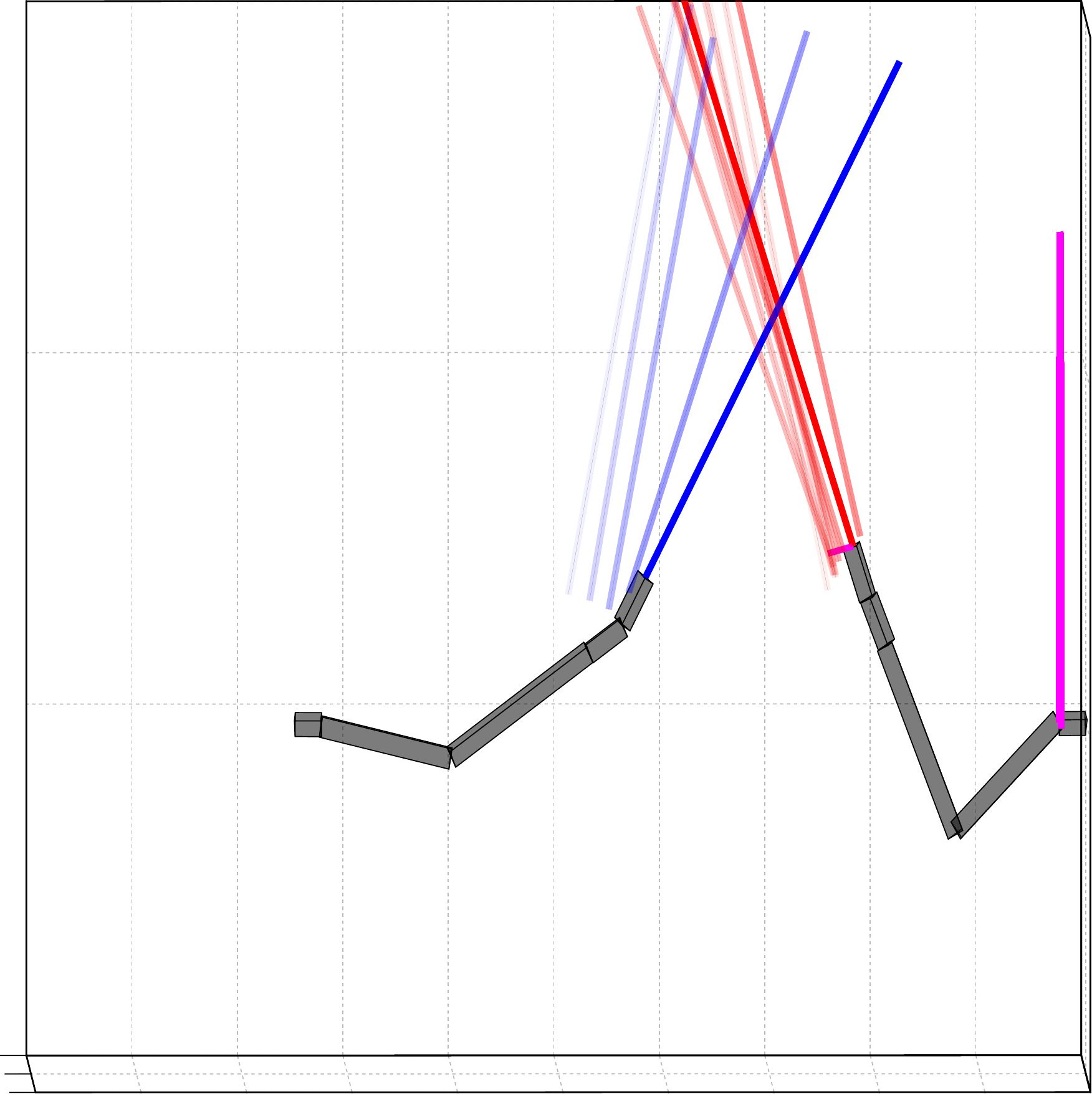} \\[3pt]
    $t = 0\,\sqparen{\textrm{ms}}$ &
    $t = 40\,\sqparen{\textrm{ms}}$ &
    $t = 80\,\sqparen{\textrm{ms}}$ &
    $t = 120\,\sqparen{\textrm{ms}}$ &
    $t = 160\,\sqparen{\textrm{ms}}$ \\
  \end{tabular}
  \caption{%
    Snapshots of the simulation in Section~\ref{subsec:example2}.
    The attacker's sword is shown in blue, the defender's sword in red, and the set $S_\text{target}$ is shown in purple.
    Each column corresponds to a different time, the top row from an inclined point of view, the bottom row from a lateral point of view.
  }
  \label{fig:ex2_trajectory}
\end{figure*}

This example shows the behavior of our algorithm when the assumption that $u_a(t)$ is constant does not hold.
The initial condition is identical to the previous example, and the opponent's trajectory is defined by:
\begin{equation}
  u_a(t) =
  \begin{cases}
    \sqmatrixs{-68 & -200 & -400 & -355 & 280 & 0}^T, & \text{if}\ t < 85\, [\textrm{ms}],\\
    \sqmatrixs{34 & -309 & -382 & 244 & -329 & 0}^T, & \text{if}\ t \geq 85\, [\textrm{ms}].
  \end{cases}
\end{equation}
Figure~\ref{fig:ex2_trajectory} shows the optimal trajectory when we do explicitly control for observation.
In this case, the robot successfully forces blade contact at time $t_f = 160\, \sqparen{\textrm{ms}}$.

\subsection{Effect of Active Perception in the Control Loop}

As described in the introduction, our motivation to implement this particular control loop was to show that it is possible to extract all the necessary 3-D information from a single camera if the data is analyzed in the context of a dynamical model, and if we can improve the capture of data at the same time we defend from the attacker.

The first objective, extraction of information using a dynamical model, is achieved by our observer (explained in Section~\ref{sec:observer}), while the second objective, dynamic improvement of the data capture by the camera, is achieved by adding the function $\phi_3$, defined in equation~\eqref{eq:phi3}, to the objective function of our optimal control problem, defined in equation~\eqref{eq:rhc}.

\begin{figure}
  \centering
  \begin{tabular}{cc}
    \fbox{%
      \includegraphics[width=.42\linewidth,trim=30 50 300 30,clip]{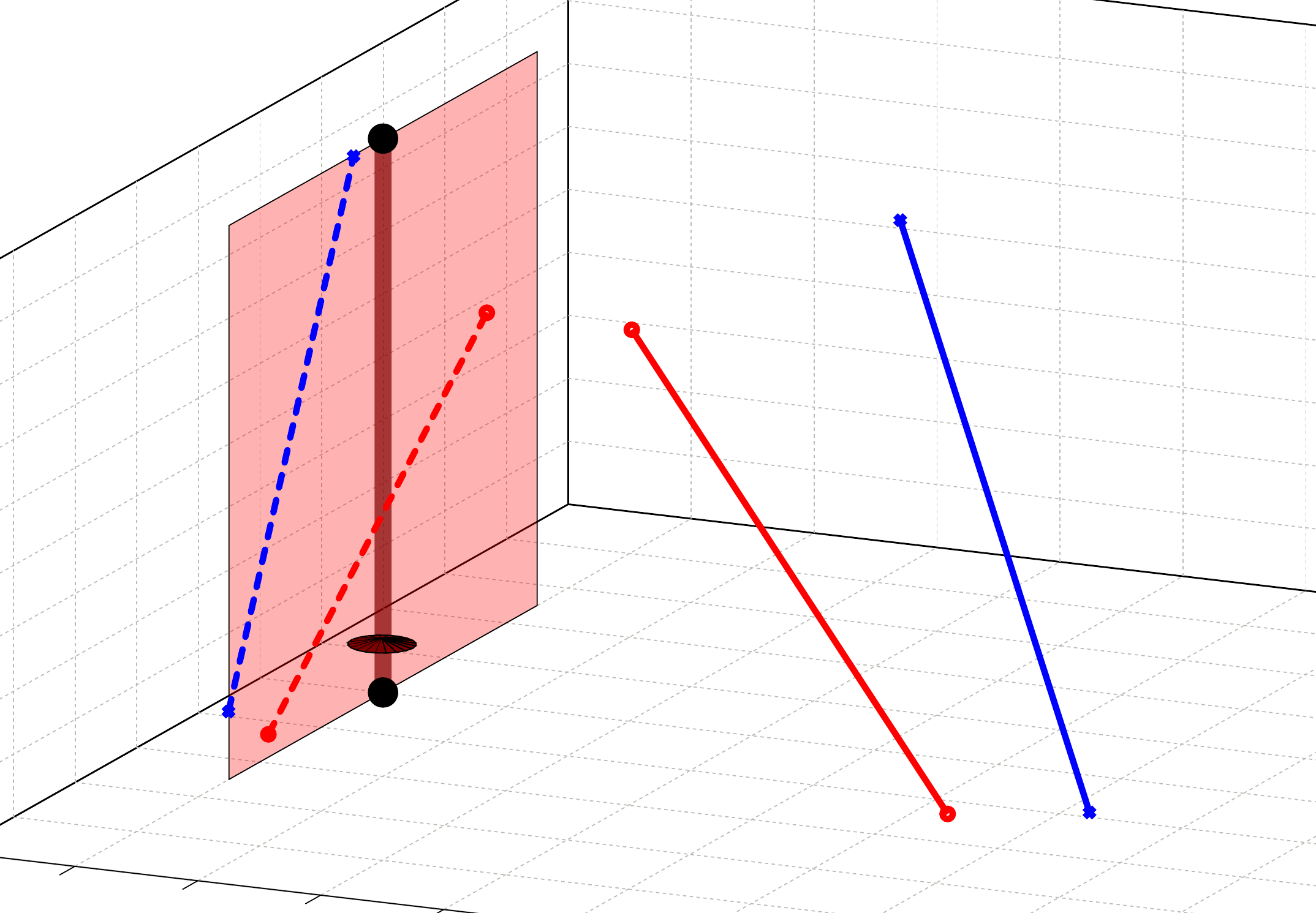}%
    } &
    \fbox{%
      \includegraphics[width=.42\linewidth,trim=30 50 300 30,clip]{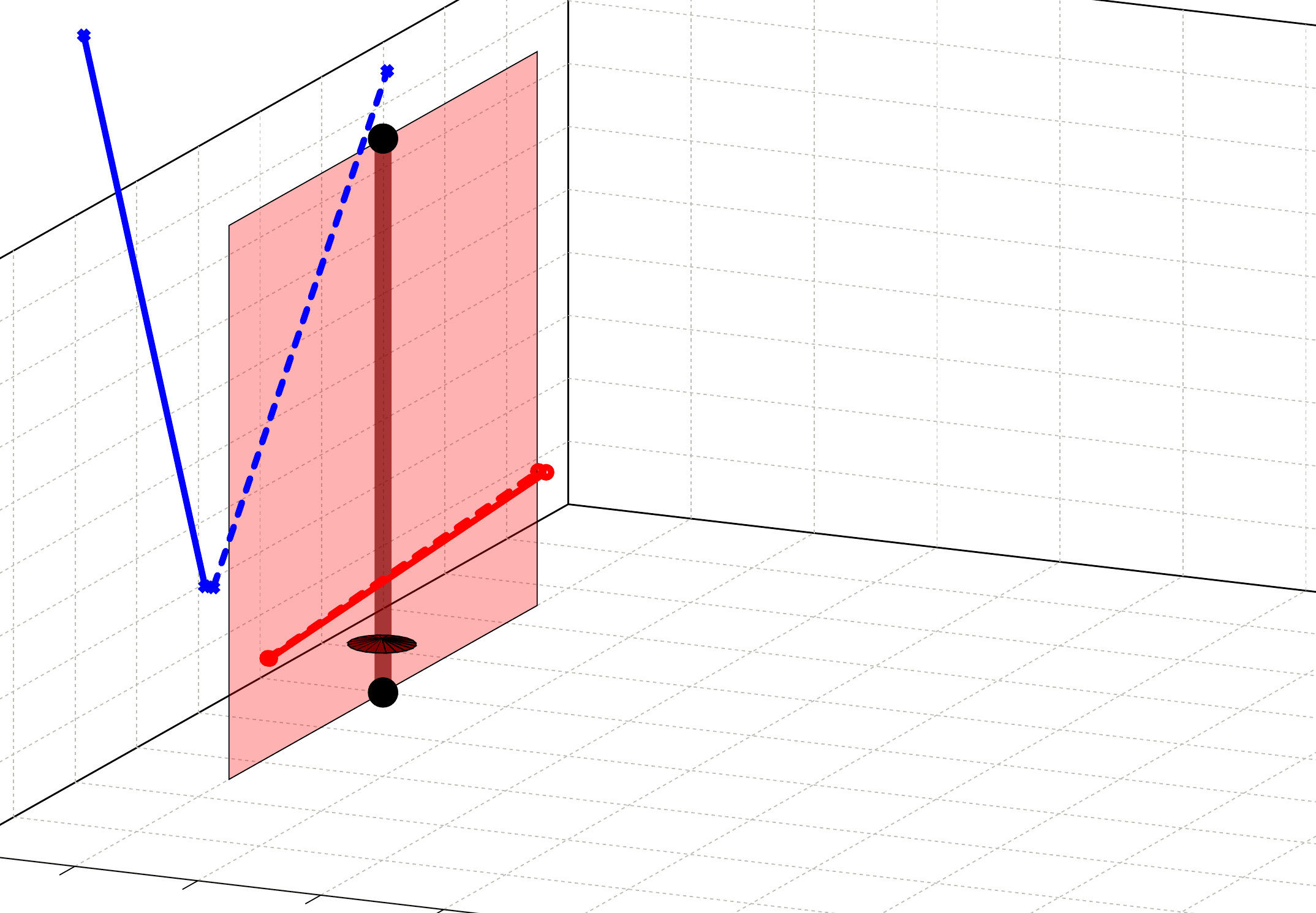}%
    } \\[5pt]
    \fbox{%
      \includegraphics[width=.42\linewidth,trim=30 50 300 30,clip]{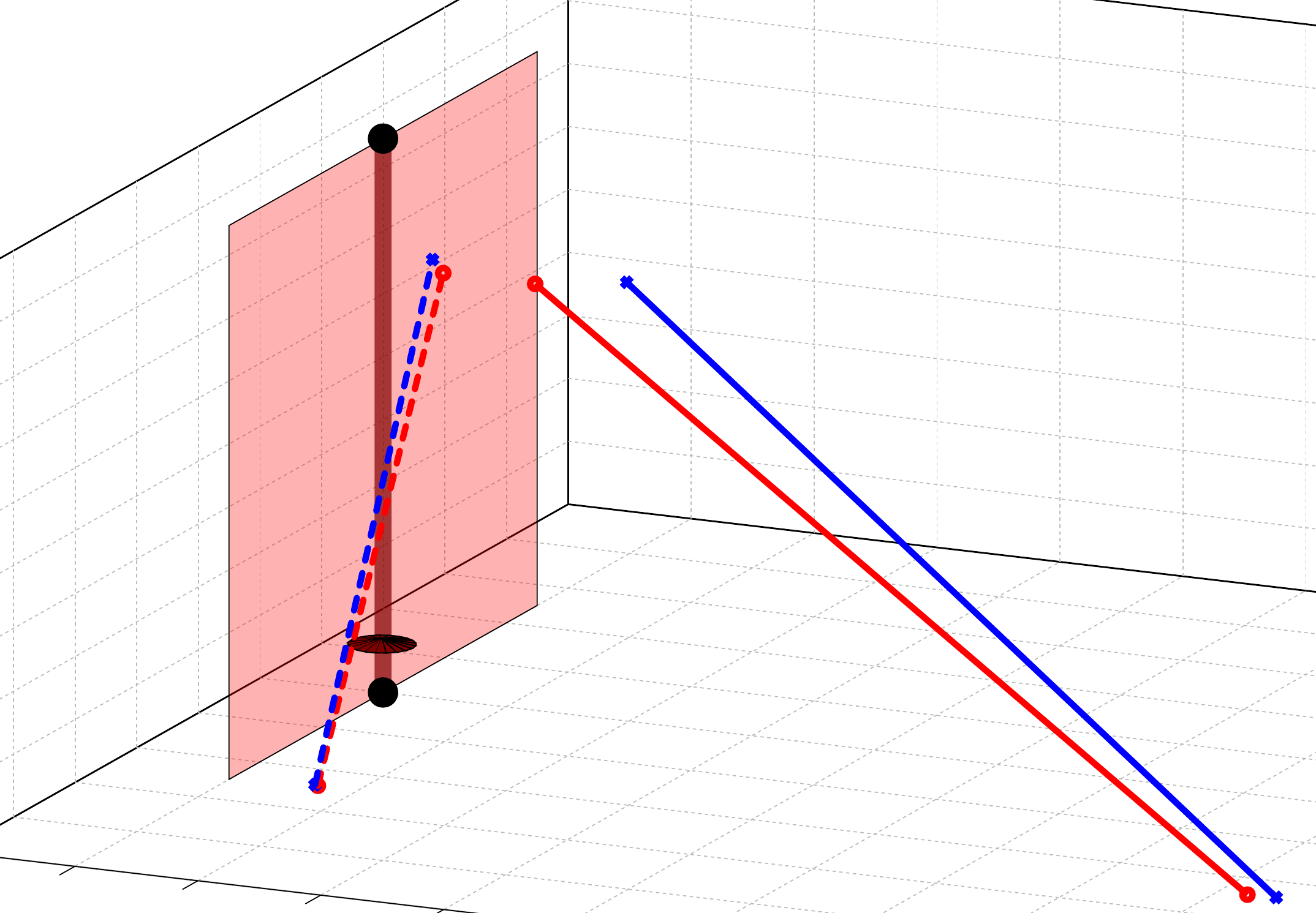}%
    } &
    \fbox{%
      \includegraphics[width=.42\linewidth,trim=30 50 300 30,clip]{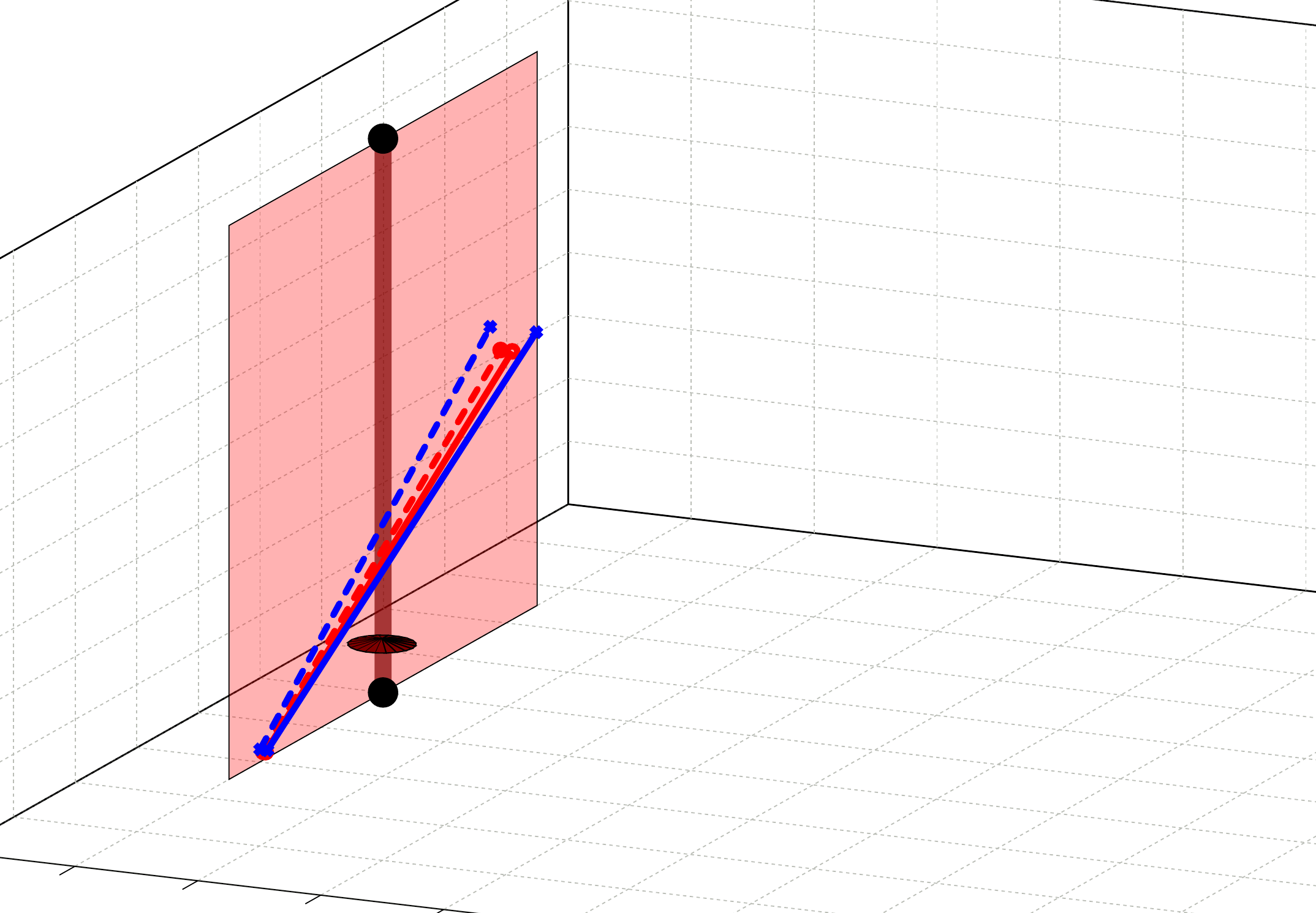}%
    } \\[3pt]
    $t = 72\,\sqparen{\textrm{ms}}$ &
    $t = 145\,\sqparen{\textrm{ms}}$ \\
  \end{tabular}
  \caption{%
    Attacker's sword (solid red), observer's estimation of its position (solid blue), both in $\paren{\axis{x}_d,\axis{y}_d,\axis{z}_d}$ coordinates, and their projections on the blocking plane (dashed red and blue, respectively), for two instants of time.
    The top row shows a simulation where $\gamma_3 = 0$, and the bottom row shows a simulation under the same conditions, but now with $\gamma_3 = 10^{-4}$.
  }
  \label{fig:ex1_observer}
\end{figure}

We evaluated the relevance of $\phi_3$ in the closed-loop performance by simulating the configuration described in Section~\ref{subsec:example1} with $\gamma_3 = 0$ and with $\gamma_3 > 0$ (using the value in Table~\ref{tab:params}).
Figure~\ref{fig:ex1_observer} shows the attacker's sword next to our estimation for both cases, $\gamma_3 = 0$ and $\gamma_3 > 0$.
For the configuration described in Section~\ref{subsec:example2} the camera quickly looses track of the attacker's sword, so the optimization problems did not produce meaningful results for $\gamma_3 = 0$.

\begin{figure}
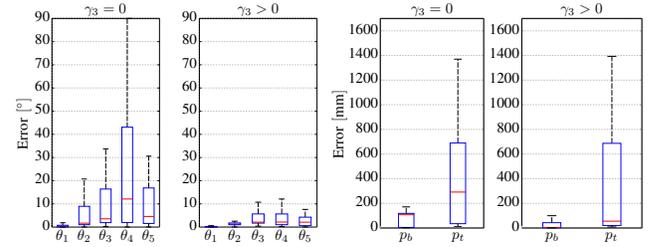

  \centering
  \mbox{}%
  \hfill%
  \resizebox{!}{.36\linewidth}{\input{figures/obsError_joints.pdftex_t}}%
  \hfill%
  \resizebox{!}{.36\linewidth}{\input{figures/obsError_points.pdftex_t}}%
  \hfill%
  \mbox{}
  \caption{%
    Distribution of the estimation error when $\gamma_3 = 0$ vs.\ $\gamma_3 > 0$ for two measurements: attacker's joint angles (left), and position of the attacker's sword tip and base (right).
    Note that the attacker's angle $\theta_6$ is omitted since it rotates the sword, hence it is irrelevant for control purposes.
  }
  \label{fig:obs_error}
\end{figure}

We also quantified the effect of explicitly controlling for observation by looking at the error the observer makes when $\gamma_3 = 0$ and $\gamma_3 > 0$.
Figure~\ref{fig:obs_error} shows the estimation error empirical distribution in each of the attacker's joint angles, and in the positions of tip and base of the attacker's sword.
The distribution was estimated using 20 different simulations, each with a different attacker input $u_a$ chosen such that the attacker's sword successfully reaches a point in $S_{\text{target}}$ provided the defender does not move.
Although the worst case errors in Figure~\ref{fig:obs_error} are comparable, particularly in the location of the tip of the sword, the median error is much lower when $\gamma_3 > 0$.
In fact, the worst case arises when the camera loses sight of the attacker's sword, but this situation occurs significantly less when $\gamma_3 > 0$.

\section{Conclusion}
\label{sec:conclusion}

We have presented a method for simultaneous sensing and control of a robotic fencing arm with a single camera attached to it.
We extracted the useful information from the images captured by the camera using a Receding Horizon Estimation algorithm, and we controlled the robotic arm to defend against an attacker using a Receding Horizon Control algorithm.
We included in our Receding Horizon Control algorithm a term to maintain the attacker's sword within the field of view or our camera, thus actively improving the data captured, and by consequence, the quality of our estimations.

Our focus in the near future includes the real-time implementation of our closed-loop system using a real robotic arm, and the use of an heterogeneous collection of sensor, rather than a single camera.

\section{Acknowledgments}

We are grateful to Tathagata Banerjee, who helped debugging our code.

\bibliographystyle{IEEEtran}
\bibliography{refs}

\end{document}